\documentclass[reqno,10pt]{amsart}
\usepackage{amsmath,amssymb,mathrsfs,amsthm,amsfonts}
\usepackage[usenames,dvipsnames]{xcolor}
\usepackage{comment}
\usepackage{hyperref}
\usepackage{graphicx}
\usepackage{enumerate}
\usepackage{stmaryrd}
\usepackage[all]{xy}
\usepackage{bbm}
\usepackage{multicol}
\usepackage{hyperref}
\usepackage[capitalise]{cleveref}
\usepackage{caption}
\usepackage{subcaption}
\hypersetup{
     colorlinks=true,
     linkcolor=blue,
     filecolor=black,
     citecolor = blue,      
     urlcolor=blue,
     }
\usepackage{fullpage}

\usepackage{tikz}
\usetikzlibrary{decorations.markings}
\usetikzlibrary{arrows}
\usetikzlibrary{calc}
\usepackage{thmtools}

\newtheorem{theorem}{Theorem}[section]
\newtheorem{lemma}[theorem]{Lemma}

\newtheorem{proposition}[theorem]{Proposition}
\newtheorem{corollary}[theorem]{Corollary}
\theoremstyle{definition}
\newtheorem{definition}[theorem]{Definition}
\newtheorem{assumption}[theorem]{Assumption}
\newtheorem{remark}[theorem]{Remark}

\numberwithin{equation}{section}
\makeatletter
\AddToHook{cmd/appendix/before}{\def\cref@section@alias{appendix}\def\cref@subsection@alias{appendix}}
\makeatother

\graphicspath{{Images/}}

\theoremstyle{plain}
\newcommand\thankssymb[1]{*}
\begin{document}
\title{Stationary Measure of the Open KPZ Equation \\ through the Enaud-Derrida Representation}
\author[Himwich]{Zoe Himwich}
\address{Department of Mathematics, Columbia University, 2990 Broadway, New York, NY 10027, USA}
\email{himwich@math.columbia.edu}
\date{\today}
\begin{abstract} 
    Recent works of Barraquand and Le Doussal \cite{BLD22} and Bryc, Kuznetsov, Wang, and Wesolowski \cite{BKWW23} described the open KPZ stationary measure as the law of the sum of two independent processes: an ordinary Brownian motion and a process whose law is absolutely continuous with respect to Wiener measure, with an explicit Radon--Nikodym derivative. Subsequent work of Barraquand and Le Doussal \cite{BD23} used the Enaud--Derrida \cite{ED04} representation of the DEHP algebra to formulate the open ASEP stationary measure in terms of the sum under a Radon--Nikodym reweighting of the law of a two-dimensional simple random walk. They show that this Radon-Nikodym derivative converges pointwise to the Radon-Nikodym derivative that characterizes the open KPZ stationary measure. This article proves that the corresponding sequence of measures converges weakly to the open KPZ stationary measure. This provides an alternative proof of the probabilistic formulation of the open KPZ stationary measure, which avoids dealing explicitly with finite-dimensional distributions. We also provide the first construction of the measure on intervals of a general length and for the full range of parameters in the fan region $(u+v>0)$.
\end{abstract}
\maketitle
%\noindent\makebox[\linewidth]{\rule{\textwidth}{.08pt}}
% \setcounter{tocdepth}{1}{
%   \hypersetup{linkcolor=black}
%   \tableofcontents
% }

\section{Introduction}

The KPZ equation is a nonlinear stochastic partial differential equation (SPDE) originally introduced by Kardar, Parisi, and Zhang~\cite{KPZ86} to describe the behavior of random interfaces under relaxation and lateral growth. By the open KPZ equation, we refer to the same SPDE with spatial coordinate restricted to an interval of finite length, and with fixed boundary conditions. For $t\geq 0, x\in[0,L]$, the equation takes the form 
\begin{align}\label{e:kpz}
    \partial_{t}H_{u,v}(t,x) = \partial_{xx}H_{u,v}(t,x) + \left(\partial_{x}H_{u,v}(t,x)\right)^{2} + \sqrt{2} \xi(t,x),
\end{align} where $\xi(t,x)$ is a space-time white noise.
The boundary conditions are given by $\partial_{x}H_{u,v}(t,x)|_{x=0}=u,$ and  $\partial_{x}H_{u,v}(t,x)|_{x=L}=-v$ for  $u,v\in\mathbb{R}$. We sometimes refer to these as Neumann boundary conditions. The equation above needs some additional explanation in order to be well-posed. In particular, we need to make sense of the nonlinear term $(\partial_{x}H_{u,v}(t,x))^{2}$. A work of Mueller~\cite{M91} demonstrates that any mild solution of the multiplicative stochastic heat equation, 
\begin{align*}
    \partial_{t}Z(t,x) = \partial_{x}^{2}Z(t,x)+\sqrt{2}Z(t,x)\xi(t,x),
\end{align*} which is almost surely positive at time zero $(Z(0,x)>0)$, will then remain positive for all $x$ and $t>0$. Consequently, the ``Hopf-Cole'' transformation $H(t,x)=\log{\left(Z(t,x)\right)}$ gives a well-defined solution to the KPZ equation (the standard equation without open boundaries). Solutions to the open KPZ equation are defined via the same Hopf-Cole transform which we use in the standard KPZ setting, with the added condition that the solution to the multiplicative stochastic heat equation on $[0,L]$, now parameterized in terms of boundary conditions $u,v\in\mathbb{R}$, must also satisfy $\partial_{x}Z_{u,v}(t,x)|_{x=0} = \left(u-\frac{1}{2}\right)Z_{u,v}(t,0),$ and $\partial_{x}Z_{u,v}(t,x)|_{x=L} = -\left(v-\frac{1}{2}\right)Z_{u,v}(t,L).$ A well-defined notion of what it means to be a solution of the open KPZ equation appears in Corwin and Shen~\cite{CS18}. 

For $L>0$, we use the notation $C_{0}[0,L]$ to denote continuous functions on
$[0,L]$, started at $0$.
Both spaces are equipped with the uniform topology.
\begin{definition}\label{d:kpzstationary} A probability measure $\mu$ on $C_{0}[0,L]$ is a stationary measure of the open KPZ increment process with Neumann boundary parameters $u,v\in\mathbb{R}$ if, whenever the initial data $\{H_{u,v}(0,x)-H_{u,v}(0,0)\}_{x\in[0,L]}$ has law $\mu$, then for all $t\geq 0$, $\{H_{u,v}(t,x)-H_{u,v}(t,0)\}_{x\in[0,L]}$ also has law $\mu$. When discussing stationary measures, we will often suppress the dependence on time and write the process with law $\mu$ as $\{H_{u,v}(x)-H_{u,v}(0)\}_{x\in [0,L]}$. 
\end{definition}
Corwin and Knizel~\cite{CK21} proved the existence of a stationary measure for the open KPZ equation for all $u,v\in\mathbb{R}$. They also characterized the measure in terms of a multipoint Laplace transform under the condition $u+v\geq 0$. Later works by Barraquand and Le Doussal~\cite{BLD22} and Bryc, Kuznetsov, Wang, and Wesolowski~\cite{BKWW23} gave a probabilistic description of this measure in terms of a stochastic process. They showed that the stationary measure of the open KPZ increment process can be described (for $L=1$ and $u+v>0$ and $u,v>-1$) in terms of a Radon-Nikodym transform applied to the joint measure of a pair of Brownian motions. \cref{d:kpzrn} gives a precise characterization of the appropriate measure, in the slightly more general setting (for any $L$, and $u+v>0$) needed to state our main theorem.

We use $\mathbb{W}_{L}$ to denote the law of a pair of independent Brownian motions satisfying $\mathbb{E}[B_{i}(s)B_{j}(t)]=2^{-1}\delta_{ij}\min\{s,t\}$ for $i,j\in\{1,2\}$. We likewise use the notation $\mathbb{L}$ to denote the Lebesgue measure on $\mathbb{R}$, as well as $\mathbb{P}_{L}:=\mathbb{L}\otimes\mathbb{W}_{L}$ to denote the product of these two measures. In the next definition, and elsewhere in the article, we will use the convention that when $\mu$ is an infinite measure, $\mathbb{E}_{\mu}[\cdot]$ refers to integrating with respect to $\mu$.

\begin{definition}\label{d:kpzrn} For $L>0$ and $u,v\in\mathbb{R}$, we define a function $\Phi_{L;u,v}(x,g,h):\mathbb{R}\times C_{0}[0,L]\times C_{0}[0,L]\to \mathbb{R}$ and a partition function $\mathcal{Z}_{L;u,v},$ 
    \begin{align*}
    \Phi_{L;u,v}(x,g,h):= \exp{\left(-2(u+v)x -2vg(L) - e^{-2x}\int_{0}^{L}e^{-2g(t)}dt\right)}, & &  \mathcal{Z}_{L;u,v} := \mathbb{E}_{\mathbb{P}_{L}}\left[\Phi_{L;u,v}(x,g,h)\right].
    \end{align*} In \cref{c:finitemeas}, we prove that $\mathcal{Z}_{L;u,v}<\infty$ for all $u+v>0$ (it follows immediately from the definition that $\mathcal{Z}_{L;u,v}>0$). For the same range of parameters, we define the probability measure $\mathbb{Q}_{L;u,v}$ through the Radon-Nikodym derivative $$\frac{d\mathbb{Q}_{L;u,v}}{d\mathbb{P}_{L}} =\mathcal{Z}_{L;u,v}^{-1}\Phi_{L;u,v}(x,g,h).$$
    \end{definition}
\begin{remark}\label{r:thirdcoordinate}
    To recover the expression in~\cite[Equation (40)]{BD23}, set $U(t)=g(t)+x$. It is advantageous for the calculations in this paper to separate the starting value $x$ from the function $g(t)$. The weight $\Phi_{L;u,v}(x,g,h)$ does not depend on its third  argument. We also note that under $\mathbb{W}_L$, the two Brownian coordinates are independent and therefore if $(X,G_{L;u,v},B_{L})\sim\mathbb{Q}_{L;u,v},$ then $B_{L}$ is independent of $(X,G_{L;u,v})$. In particular, $B_{L}$ is a Brownian motion with diffusion coefficient $1/\sqrt{2}$. We retain $h$ in the expression because both path coordinates enter the representation of the stationary profile in
\cref{t:main1}.
\end{remark}
    We can now state the main result. 
\begin{theorem}\label{t:main1} Let $(X,G_{L;u,v}, B_{L})\sim \mathbb{Q}_{L;u,v}$. Then the law of $\{B_{L}(x)+G_{L;u,v}(x)\}_{x\in[0,L]}$ is the stationary measure of the open KPZ increment process on the interval $[0,L]$ with Neumann boundary parameters $u,v\in\mathbb{R}$ such that $u+v>0$. 
\end{theorem}    
That there is a unique such measure is due to the work of Knizel and Matetski~\cite{KM22} and Parekh~\cite{P22}. That this is a stationary measure for the open KPZ increment process when $L=1$ and $u+v>0$, $u,v>-1$ is known from \cite{BLD22,BKWW23}. These articles predicted that the same would hold in the general case that $u+v>0$. We prove \cref{t:main1} for general $L$ and all $u+v>0$, and  provide a rigorous proof of a new derivation of the stationary measure through the Enaud-Derrida representation of the DEHP algebra~\cite{ED04} (see \cref{s:MPAED} for background on the open ASEP, the DEHP algebra, and the Enaud-Derrida representation). 
    
    Barraquand and Le Doussal~\cite{BD23} use the Enaud-Derrida representation to describe the stationary measure of the open ASEP increment process (\cref{d:asepstationary}) as the sum of reweighted random walk processes. They demonstrate that, under appropriate scaling, the functional which reweights the processes converges pointwise to $\Phi_{L;u,v}(x,g,h)$ (\cref{t:pointwise}). We build on their calculation to show weak convergence of measures. This approach has the advantage of circumventing the several steps involved in the original construction of the open KPZ stationary measure. As noted above, it also provides a probabilistic construction of the open KPZ stationary measure in a wider range of parameters. This article also provides an alternative proof that the partition function in \cref{d:kpzrn} is finite, a fact previously proved in \cite{BKWW23} for $L=1$. We provide an alternative argument (for general $L$), following a paper of Yor~\cite{Y92}. 
\subsection{Outline and Discussion}\label{ss:outline} In \cref{s:tech}, we introduce the open ASEP, the matrix product ansatz, and the Enaud-Derrida representation. We discuss works of Derrida, Enaud, and Lebowitz~\cite{DEL} and Barraquand and Le Doussal~\cite{BD23} which use the Enaud-Derrida representation to describe the stationary measure of the open ASEP increment process as the sum of reweighted random walks. We establish the rescaling that is needed to obtain the stationary measure of the open KPZ increment process from their formula. In the process, we define a sequence of measures $\mathbb{Q}^{(N)}_{L;u,v}$ (\cref{d:qqn}) which are related to the rescaled stationary measure of open ASEP in the same way as $\mathbb{Q}_{L;u,v}$ is related to the stationary measure of open KPZ. We state the main technical result (\cref{t:main}), that $\mathbb{Q}^{(N)}_{L;u,v}$ converge weakly to $\mathbb{Q}_{L;u,v}$. Finally, we prove \cref{t:main1} from \cref{t:main}. 

\begin{remark}
    Before explaining the structure of the subsequent sections, we will briefly comment on the method used to prove \cref{t:main} and the main challenges. In this remark, we will speak about the random variables and measures of interest in general terms; full definitions will follow in \cref{s:tech}. 
We have a collection of random variables $Y^{(N)}$ with law $\mathbb{Q}^{(N)}_{L;u,v}$. This law is obtained by reweighting a reference measure by a function $\Phi_{L;u,v}^{(N)}$ and then normalizing. In particular, we let
\begin{align*}
 \mathbb{L}^{(N)} := \frac{1}{\sqrt N}\sum_{j\geq 1}\delta_{j/\sqrt N - \log{\sqrt{N}}},
\end{align*} and let $\mathbb{W}_L^{(N)}$ denote the law of a $\lfloor LN\rfloor$-step two-dimensional simple symmetric random walk started at $(0,0)$, scaled in its two coordinate values by $N^{-1/2}$ and linearly interpolated between time values $t_i:=iL\lfloor LN\rfloor^{-1}$ for $i\in \llbracket 0,\lfloor LN\rfloor\rrbracket$. We define a reference measure $\mathbb{P}_{L}^{(N)}:=\mathbb{L}^{(N)}\otimes \mathbb{W}_L^{(N)},$ from which we obtain $\mathbb{Q}^{(N)}_{L;u,v}$ through reweighting by $\Phi_{L;u,v}^{(N)}$ and normalizing. Recall that we use $\mathbb{L}$ to denote Lebesgue measure on $\mathbb{R}$, and $\mathbb{W}_L$ to denote the law of a pair of independent Brownian motions $(B_1,B_2)$ satisfying $\mathbb{E}[B_i(s)B_j(t)]=2^{-1}\delta_{ij}\min\{s,t\},$ and the definition $\mathbb{P}_{L}:=\mathbb{L}\otimes\mathbb{W}_L.$ The lattice measures $\mathbb{L}^{(N)}$ converge vaguely to $\mathbb{L}$, while the random-walk path measures $\mathbb{W}_L^{(N)}$ converge weakly to $\mathbb{W}_L$. Neither $\mathbb{P}_{L}^{(N)}$ nor $\mathbb{P}_{L}$ is a probability measure. We define
\begin{align*}
\mathcal{Z}^{(N)}_{L;u,v} := \mathbb{E}_{\mathbb{P}_{L}^{(N)}}\left[\Phi_{L;u,v}^{(N)}\right].
\end{align*}
The candidate limiting Radon-Nikodym derivative reweighting ($\Phi_{L;u,v}$ and its normalizing constant $\mathcal{Z}_{L;u,v}$) are given in \cref{d:kpzrn}. If $Y^{(N)}$ has law $\mathbb{Q}_{L;u,v}^{(N)}$ and $Y$ has law $\mathbb{Q}_{L;u,v}$, then a typical criterion to establish weak convergence is that for
every bounded continuous function $F$, 
\begin{multline}\label{e:heur}
\lim_{N\to\infty}\mathbb{E}[F(Y^{(N)})] = \lim_{N\to\infty}(\mathcal{Z}^{(N)}_{L;u,v})^{-1}\mathbb{E}_{\mathbb{P}_{L}^{(N)}}\left[ F \cdot \Phi_{L;u,v}^{(N)}\right] = (\mathcal{Z}_{L;u,v})^{-1}\mathbb{E}_{\mathbb{P}_{L}}\left[ F \cdot \Phi_{L;u,v}\right] = \mathbb{E}[F(Y)].
\end{multline} 
To prove \eqref{e:heur}, we establish the following:
\begin{enumerate}
    \item For every bounded Lipschitz function $F$,
    \begin{align*}
        \lim_{N\to \infty}\left|\mathbb{E}_{\mathbb{P}_{L}^{(N)}}\left[F\cdot\Phi_{L;u,v}^{(N)}\right]-\mathbb{E}_{\mathbb{P}_{L}}\left[F\cdot\Phi_{L;u,v}\right]\right|=0,
    \end{align*}
    which we prove in two parts:
    \begin{enumerate}
        \item $\lim_{N\to\infty}\big|
        N^{-1/2}\sum_{x\in\mathbb Z^{(N)}}\mathbb{E}_{\mathbb W_L}\left[F(x,\cdot,\cdot)\Phi_{L;u,v}(x,\cdot,\cdot)\right]-\mathbb{E}_{\mathbb P_L}\left[F\cdot\Phi_{L;u,v}\right]\big| = 0;$
        \item $\lim_{N\to \infty} N^{-1/2}\sum_{x\in\mathbb Z^{(N)}} \big|\mathbb{E}_{\mathbb W_L^{(N)}}\left[F(x,\cdot,\cdot)\Phi_{L;u,v}^{(N)}(x,\cdot,\cdot)\right]-\mathbb{E}_{\mathbb W_L}\left[F(x,\cdot,\cdot)\Phi_{L;u,v}(x,\cdot,\cdot)\right]\big| =0.$
    \end{enumerate}
    \item $\mathcal Z_{L;u,v}$ is finite.
\end{enumerate}
We note that combining $(1)$ (with $F\equiv \mathbf{1}$) and $(2)$ allows us to conclude $\mathcal{Z}^{(N)}_{L;u,v}\to\mathcal{Z}_{L;u,v}$. 
We rely on two key bounds, \cref{l:finitebounds} and \cref{l:limitbounds}, at several points in the proof of $(1a-b)$. These two lemmas are the main technical challenge of the paper. The additional complexity added from the argument of \cite{BD23} are exactly these more sophisticated controls on the Radon-Nikodym derivatives. 

\cref{l:finitebounds} can also be applied to show that the $\mathcal{Z}^{(N)}_{L;u,v}$ are bounded by a constant, which, along with $(1)$, would imply $(2)$. \cref{l:limitbounds} also directly implies $(2)$. There are several ways one might proceed to prove a result of the form $(1b)$, we use the intuitive and direct approach of a KMT embedding, after separating out the dependence on the random initial point from the random walk and Brownian paths. 
\end{remark}

 In \cref{s:finite}, we prove \cref{l:finitebounds} and \cref{l:limitbounds}. We also apply \cref{l:limitbounds} to show that $\mathbb{Q}_{L;u,v}$ is a probability measure, which involves showing that $\mathcal{Z}_{L;u,v}$ is finite. In \cref{s:conv}, we establish a bound on the distance between $\mathbb{E}_{\mathbb{W}^{(N)}_{L}}[\Phi^{(N)}_{L;u,v}]$ and $\mathbb{E}_{\mathbb{W}_{L}}[\Phi_{L;u,v}]$ as a function of $x$ (analogous to $(1b)$ above). In \cref{s:weakconvergencefinal} we use the bound from \cref{s:conv} to prove \cref{t:main}. 
\cref{s:rws} contains the proof of several lemmas which are used in \cref{s:finite}, and \cref{a:pointwise} contains the proof of pointwise convergence of the Radon-Nikodym derivatives defined in \cite{BD23} to $\Phi_{L;u,v}(x,g,h)$. 

\subsection{Acknowledgements} The author thanks Ivan Corwin for suggesting this question and for comments on the draft. The author particularly thanks Zongrui Yang for many helpful discussions and comments. The author thanks Shalin Parekh for additional discussions. The author was supported by the Fernholz Foundation's Summer Minerva Fellows Program, as well as Ivan Corwin's grant, NSF DMS-1811143. 

\section{Background and Main Technical Result}\label{s:tech}

In this section, we build up the definitions and notation which we need to state the main technical result which goes into the proof of \cref{t:main1} (\cref{t:main}). We conclude this section by giving the proof of \cref{t:main1} from \cref{t:main}, and further sections of the paper will be devoted to the proof of \cref{t:main}. We begin by defining some notation.

\subsection{Notation}\label{ss:Notation}

For $a,b\in\mathbb{R}$ with $a\leq b$, we write $\llbracket a,b\rrbracket := [a,b]\cap\mathbb{Z},$ or, for $b\geq 1$, $\llbracket b\rrbracket:=[1,b]\cap\mathbb{Z}$ We use $N\in\mathbb{Z}_{>0}$ to denote the weak-asymmetry and inverse-macroscopic-mesh parameter. For each $N$ such that $\lfloor LN\rfloor\geq1$, we define $t_i:=iL\lfloor LN\rfloor^{-1},$ for $i\in\llbracket0,\lfloor LN\rfloor\rrbracket.$ The $N$-th finite open ASEP system will have $\lfloor LN\rfloor$ sites, and the associated two-dimensional random walk will have $\lfloor LN\rfloor$ steps. The rescaled path will have values scaled by $N^{-1/2}$ and will be linearly interpolated between $\{t_i\}_{i=0}^{\lfloor LN\rfloor}$.

For $n\in\mathbb{Z}_{\geq 0}$ and $q\in[0,1)$, we define the $q$-bracket $[n]_q:=(1-q^n)/(1-q).$

When discussing the matrix product ansatz and the Enaud--Derrida representation, we use bra--ket notation, with $\langle\cdot|$ denoting a row vector and $|\cdot\rangle$ denoting a column vector.

We write $C[0,L]$ for the space of real-valued continuous functions on $[0,L]$, and define the notation $C_{0}[0,L]:=\left\{f\in C[0,L]:f(0)=0\right\}$. Both spaces are equipped with the uniform topology. For $T>0$, we define $C([0,T],C[0,L])$ to be the space of continuous functions from $[0,T]$ into $C[0,L]$, equipped with the uniform topology. We similarly define $D([0,T],C[0,L])$ to be the space of cadlag functions from $[0,T]$ into $C[0,L]$, equipped with the Skorokhod topology induced by the uniform metric on $C[0,L]$.

For $\vec f=(f_1,f_2)\in C_0[0,L]\times C_0[0,L]$, we define $\|\vec f\|_{\infty,2}:=\sup_{t\in[0,L]}\|\vec f(t)\|_{\mathbb R^2}.$ We equip $\mathbb R\times C_0[0,L]\times C_0[0,L]$ with the metric $d\big((x,g,h),(x',g',h')\big):=|x-x'|+\|(g-g',h-h')\|_{\infty,2},$ and the usual product topology. For a function
$F:\mathbb R\times C_0[0,L]\times C_0[0,L]\to\mathbb R$, we define $$\mathrm{Lip}(F):=\sup_{\substack{(x,g,h)\neq(x',g',h')}}\frac{|F(x,g,h)-F(x',g',h')|}{d\big((x,g,h),(x',g',h')\big)}.$$ We say that $F$ is bounded Lipschitz if $\|F\|_{\infty}+\mathrm{Lip}(F)<\infty$.

We use $\mathbb{L}$ to denote the Lebesgue measure on $\mathbb{R}$ and we use $\mathbb{W}_{L}$ to denote the law of a pair of independent Brownian motions satisfying $\mathbb{E}[B_{i}(s)B_{j}(t)]=2^{-1}\delta_{ij}\min\{s,t\}$ for $i,j\in\{1,2\}$. We use the notation $\mathbb{P}_{L}:=\mathbb{L}\otimes \mathbb{W}_{L}$ for the product measure. As before, we write $\mathbb{L}^{(N)} := \sum_{i\geq 1} N^{-\frac{1}{2}} \delta_{i/\sqrt{N}-\log{\sqrt{N}}}.$ We use the notation $\mathbb{E}_{P}[\cdot]$ to denote integrating with respect to the measure $P$. The measure $P$ is not always a probability measure (for example, $\mathbb{P}_{L}$), so this is a slight abuse of notation. 
 
\subsection{The open ASEP}
In our notation, the $N$-th open ASEP will be defined with underlying lattice $\llbracket\lfloor LN\rfloor\rrbracket$. We use $\alpha,\gamma$ to denote the rate of particle movement into and out of the leftmost lattice point from, respectively to, the left reservoir. Similarly, we use $\beta,\delta$ to denote the rate of particle movement out of, respectively into, the rightmost lattice point of the lattice to, respectively from, the right reservoir. Within the lattice, particles move right at rate $1$ and left at rate $q\in[0,1)$, under the exclusion constraint that a particle will not jump left or right if there is already a particle occupying the destination site. We often reparameterize the open ASEP in terms of Askey-Wilson parameters $A,B,C,D$, 
\begin{align*}
     A := \kappa^{+}(q,\beta,\delta), & & B := \kappa^{-}(q,\beta,\delta), & & C := \kappa^{+}(q,\alpha,\gamma), & & D := \kappa^{-}(q,\alpha,\gamma),
 \end{align*}
 \begin{align*}
     \kappa^{\pm}(q,x,y) = \frac{1}{2x}\left(1-q-x+y\pm \sqrt{(1-q-x+y)^{2}+4xy}\right).
 \end{align*}
The state space of the index $N$ open ASEP is $\{0,1\}^{\lfloor LN\rfloor}$, and at a particular time $t$, $\tau(t) = (\tau_{1}(t),...,\tau_{\lfloor LN\rfloor}(t))$ is one random element of that statespace, and $\tau_{i}(t)=1$ signifies that there is a particle in position $i$ at time $t$ and $\tau_{i}(t)=0$ when position $i$ on the lattice is empty at time $t$. The open ASEP process $(\tau(t))_{t\geq 0}$ is a Markov process defined by that state space and an infinitesimal generator $\mathcal{L}$ which acts on functions of the open ASEP state space $f:\{0,1\}^{\lfloor LN\rfloor}\to\mathbb{R}$. 
The time evolution of ASEP induces a Markov process, the height function process, $h_{N;\alpha,\beta,\gamma,\delta}(t,r) = h_{N;\alpha,\beta,\gamma,\delta}(t,0) + \sum_{i=1}^{r}(2\tau_{i}(t)-1),$ for $r\in\llbracket 0,\lfloor LN\rfloor \rrbracket$ and where $h_{N;\alpha,\beta,\gamma,\delta}(t,0)$ is the net number of particles which have entered through the left boundary at time $t$. 

\begin{definition}
\label{d:asepstationary}
A probability measure $\mu$ on the space of open ASEP increments is called a stationary measure for the open ASEP increment process if, whenever $\left\{h_{N;\alpha,\beta,\gamma,\delta}(0,r)-h_{N;\alpha,\beta,\gamma,\delta}(0,0)\right\}_{r\in\llbracket0,\lfloor LN\rfloor\rrbracket}$ has law $\mu$, then, for every $t\geq 0$, $\left\{h_{N;\alpha,\beta,\gamma,\delta}(t,r)-h_{N;\alpha,\beta,\gamma,\delta}(t,0)\right\}_{r\in\llbracket 0,\lfloor LN\rfloor\rrbracket}$
also has law $\mu$.
\end{definition}
It is convenient to describe the open ASEP phase diagram in terms of the effective left and right boundary densities $\rho_\ell$ and $\rho_r$. See the survey \cite{C22} for more information about the boundary densities. For the purposes of this article, we use only the
identities
\begin{align}\label{e:rhos}
    \rho_{\ell}=\frac{1}{1+C}, & &\rho_{r}=\frac{A}{1+A}.
\end{align}
\cite{CK21} constructed stationary measures of the open KPZ equation as weak-asymmetry limits of stationary open ASEP increment measures. We use the same weak-asymmetry and parameter scaling, but our $N$th open ASEP will have $\lfloor LN\rfloor$ sites so that the macroscopic spatial interval is $[0,L]$.
\begin{assumption}\label{d:scale} For each $N\in\mathbb{Z}_{>0}$ such that $\lfloor LN\rfloor\geq1$, we use the following scaling assumptions:
    \begin{enumerate}
        \item Weak asymmetry scaling: $q=e^{-2/\sqrt{N}}$.
        \item Askey-Wilson parameter scaling: 
        \begin{align*}
            A = q^{v}, & & B = -q, & & C = q^{u}, & & D= -q.
        \end{align*}
        \item Height function scaling: For $i\in \llbracket 0,\lfloor LN\rfloor\rrbracket$ we define the rescaled stationary increment \begin{align}\label{e:rescaledstationary} h^{(N)}_{L;u,v}(t_{i}) - h^{(N)}_{L;u,v}(0) := N^{-1/2}\left( h_{N;\alpha,\beta,\gamma,\delta}(i)-h_{N;\alpha,\beta,\gamma,\delta}(0)\right).
        \end{align} We extend the domain of $h_{L;u,v}^{(N)}$ to $[0,L]$ by linear interpolation.
    \end{enumerate}
    Conditions $(1)$ and $(2)$ imply triple point asymptotics 
        \begin{align*}
            \rho_{\ell} = \frac{1}{2}+\frac{u}{2}N^{-\frac{1}{2}} + o(N^{-1}),  & & \rho_{r} = \frac{1}{2}-\frac{v}{2}N^{-\frac{1}{2}} + o(N^{-1}).
        \end{align*}
\end{assumption}
\begin{remark}
    The constraint that $B=D=-q$ is equivalent to imposing Liggett's condition,
    \begin{align*}
        \alpha+\gamma/q = 1, & & \beta + \delta/q = 1.
    \end{align*}
\end{remark}

\subsubsection{Matrix Product Ansatz Techniques and the Enaud-Derrida Representation}\label{s:MPAED}
A work of Derrida, Evans, Hakim and Pasquier~\cite{DEHP93} introduced a method, the matrix product ansatz, for obtaining the stationary distribution for the open ASEP. The idea behind this technique is to put forth the ansatz that there exists an algebra (the ``DEHP'' algebra, after the authors of the original paper) with a representation in terms of matrices $\mathbf{D}$ and $\mathbf{E}$ and vectors $\langle W|$ and $|V\rangle$ (possibly infinite-dimensional) such that the stationary probability of a given configuration $\tau\in \{0,1\}^{\lfloor LN\rfloor}$ is given by $$\mathbb{P}(\tau) = \frac{\langle W | \prod_{i=1}^{\lfloor LN\rfloor}\left(\mathbf{D}\tau_{i}+\mathbf{E}(1-\tau_{i})\right)|V\rangle}{\langle W|(\mathbf{D}+\mathbf{E})^{\lfloor LN\rfloor}|V\rangle}.$$ Imposing the condition that this distribution is stationary results in the relations 
\begin{align*}\begin{split}
    \mathbf{DE} - q\mathbf{ED}  & = \mathbf{D}+\mathbf{E}, \\
    (\beta \mathbf{D}-\delta \mathbf{E})|V\rangle & = | V\rangle, \\
    \langle W|(\alpha\mathbf{E}-\gamma \mathbf{D}) & = \langle W|.
    \end{split}
\end{align*} 
\begin{remark}
    We guarantee that the denominator $\langle W| (\mathbf{D}+\mathbf{E})^{\lfloor LN\rfloor}|V\rangle$ is nonzero by assuming that $\langle W|V\rangle\neq 0$ and choosing $ABCD\notin\{q^{-\ell}|\ell\in\mathbb{Z}_{\geq 0}\}$. In particular, this will always be satisfied in the fan region, where $AC<1$, $A,C\geq 0$, and $B,D\in (-1,0].$ See \cite{MS97,ER96} for further discussion of this constraint. 
\end{remark}Mathematicians and physicists have discovered several representations of the DEHP algebra. The original representation~\cite{DEHP93} was used to solve the open TASEP (the case $q=0$). Later, foundational papers by Sasamoto~\cite{Sa99} and Uchiyama, Sasamoto, and Wadati~\cite{USW} gave representations of the DEHP algebra in terms of Askey-Wilson polynomials. Enaud and Derrida
\cite{ED04} studied a different representation of the DEHP algebra, now called the Enaud-Derrida representation (a closely related representation was discovered independently by Corteel and Nunge~\cite{CN20}). The Enaud-Derrida representation is given by
\begin{align*}
    \mathbf{D}  = \left[\begin{matrix}
        [1]_{q} & [1]_{q} & 0 & 0 & 0 & ... \\ 0 & [2]_{q} & [2]_{q} & 0 & 0 & ... \\ 0 & 0 & [3]_{q} & [3]_{q} & 0 & ...
        \\ ... & ... & .. & ... & ... & ...
    \end{matrix}\right], & & \mathbf{E} = \left[\begin{matrix}
        [1]_{q} & 0 & 0 & 0 & ... \\ [2]_{q} & [2]_{q} & 0 & 0 & ... \\ 0 & [3]_{q} & [3]_{q} & 0 & ... \\ ... & ... & ... & ... & ...
    \end{matrix}\right], \\ 
    \langle W| = \sum_{n\geq 1}\left(\frac{1-\rho_{\ell}}{\rho_{\ell}}\right)^{n}\langle n|, & & |V\rangle = \sum_{n\geq 1}\left(\frac{\rho_{r}}{1-\rho_{r}}\right)^{n}[n]_{q}|n\rangle,
\end{align*} where $[\cdot]_{q}$ is defined in \cref{ss:Notation}, and $|n\rangle$ and $\langle n|$ denote basis vectors.
\cite{DEL} noticed that it is possible to write the stationary increment of the open ASEP as the sum of two random walks, $$\left\{h_{N;\alpha,\beta,\gamma,\delta}(r)-h_{N;\alpha,\beta,\gamma,\delta}(0)\right\}_{r\in\llbracket\lfloor LN\rfloor\rrbracket}= \left\{\sum_{j=1}^{r}(2\tau_{j}-1)\right\}_{r\in\llbracket\lfloor LN\rfloor\rrbracket} \overset{d}{=} \left\{n_{r}-n_{0}+m_{r}\right\}_{r\in\llbracket\lfloor LN\rfloor\rrbracket},$$ where the second equality is in distribution. \cite{BD23} used the Enaud-Derrida representation to define a measure on $\left(n_{r},m_{r}\right)_{r=1}^{\lfloor LN\rfloor}$ by reweighting the law of a two-dimensional random walk started at $(n_{0},0)$ for any $n_{0}\in\mathbb{Z}_{>0}$. Equivalently, writing $r:=n_0$ and recentering the first coordinate by subtracting its initial value, we may separate the starting point from the random-walk path. Thus, in the formulation below, $(n_i,m_i)_{i=0}^{\lfloor LN\rfloor}$ starts at $(0,0)$, while $r\in\mathbb{Z}_{>0}$ records the original starting value, so that the uncentered first coordinate at step $i$ is $n_i+r$.
\begin{definition}[Equation (25)~\cite{BD23}]\label{d:pnotrescaled} We denote the space of two-dimensional $\lfloor LN\rfloor$-step simple symmetric random walks started at $(0,0)$ by $W_{\lfloor LN\rfloor,0}$ and, when $u+v>0$, define an unnormalized weight associated to elements of $\mathbb{Z}_{>0}\times W_{\lfloor LN\rfloor,0}$. We set $w_{\lfloor LN\rfloor}(r,\vec n,\vec m)=0$ unless $n_i+r>0$for every $i\in\llbracket0,\lfloor LN\rfloor\rrbracket$.
On this support, we then define
\begin{align}\label{e:edweight}
        w_{\lfloor LN\rfloor}(r,\vec{n},\vec{m}) & := 4^{\lfloor LN\rfloor}\left(\frac{1-\rho_{\ell}}{\rho_{\ell}}\right)^{r}\left(\frac{\rho_{r}}{1-\rho_{r}}\right)^{n_{\lfloor LN\rfloor }+r}\prod_{i=0}^{\lfloor LN\rfloor}[n_{i}+r]_{q}\mathbb{P}^{\textsc{SSRW}}(\vec{n},\vec{m}),
    \end{align} and define $\widetilde{Z}_{\lfloor LN\rfloor}(q):=\sum_{r=1}^{\infty} \sum_{(\vec n,\vec m)\in W_{\lfloor LN\rfloor,0}}w_{\lfloor LN\rfloor}(r,\vec n,\vec m)$. Finally, we define a probability measure \begin{align}\label{e:rn2} \mathbb{P}^{\textsc{ED}}_{\lfloor LN\rfloor}(r,\vec n,\vec m):=(\widetilde{Z}_{\lfloor LN\rfloor}(q))^{-1}w_{\lfloor LN\rfloor}(r,\vec n,\vec m).\end{align} The notation $\rho_{r}$ and $\rho_{\ell}$ is defined in \eqref{e:rhos}, and $\mathbb{P}^{\textsc{SSRW}}$ is the law of a two-dimensional simple symmetric random walk started at $(0,0)$.
\end{definition}
\begin{remark}
The terms $[n_{i}+r]_{q}$ in \eqref{e:edweight} guarantee that paths with nonzero weight have $n_{i}+r>0$ for all $i\in\llbracket 0,\lfloor LN\rfloor\rrbracket$. \end{remark}
 \begin{theorem}[Section 3~\cite{BD23}]\label{t:bldmaintheorem} The stationary measure of the open ASEP increment process on the finite lattice of size $\lfloor LN\rfloor$ is given by the law of $\left\{n_{i}+m_{i}\right\}_{i\in\llbracket 0,\lfloor LN\rfloor\rrbracket}$, where $(r,\vec{n},\vec{m})$ is sampled from the measure in \cref{d:pnotrescaled}. 
    \end{theorem}   
    Barraquand and Le Doussal \cite{BD23} made several additional contributions beyond  proposing that this representation of the stationary measure of the open ASEP increment process could be used to obtain the stationary measure of the open KPZ increment process. They identified the rescaling of the Enaud--Derrida probability measure in \eqref{e:rn2} that is appropriate for obtaining the open KPZ stationary measure  \cite[Equation~(33)]{BD23}. They expressed the rescaled open ASEP law and the limiting open KPZ law as normalized reweightings of appropriate discrete and continuum reference measures, and showed that the corresponding recentered discrete weights converge pointwise to $\Phi_{L;u,v}$ (see \cref{t:pointwise}). 
    We will now describe their rescaling in preparation for the proof of \cref{t:main1}. 
For $i\in\llbracket0,\lfloor LN\rfloor\rrbracket,$ we set $x:=rN^{-1/2}$, $g(t_i):=n_i N^{-1/2}$, and $h(t_i):=m_i N^{-1/2}$. We extend the domain of $g$ and $h$ to $[0,L]$ by linear interpolation. We use $W_{[0,L]}^{(N)}$ to denote the set of paths of the form $(g,h)$ obtained through this rescaling and linear interpolation. We also use $\mathbb{W}_{L}^{(N)}$ to denote the pushforward of $\mathbb{P}^{\mathrm{SSRW}}$ under the rescaling and linear interpolation. This notation is motivated by the weak convergence of $\mathbb{W}_{L}^{(N)}$ to $\mathbb{W}_{L}$.
At this stage, the initial coordinate takes values in $ N^{-1/2}\mathbb{Z}_{>0}.$ 
The measure $\mathbb{L}^{(N)}$ defined in \cref{ss:Notation} corresponds to the recentered lattice measure, obtained by the change of variables $x\mapsto x-\log\sqrt N.$ 
We now describe the weak-asymmetry rescaling of \eqref{e:rn2}.
\begin{definition}
\label{d:rn}
Let $S_N := N^{-1/2}\mathbb{Z}_{>0} \times W_{[0,L]}^{(N)}.$ We define a weight function $R^{(N)}:S_N\longrightarrow [0,\infty)$ for all $(x,g,h)$ such that $g(t_{i})+x>0$ for all $i\in \llbracket 0,\lfloor LN\rfloor\rrbracket$ by
\begin{align*}
R^{(N)}(x,g,h) :=\left(\frac{1-\rho_\ell}{\rho_\ell}\right)^{\sqrt N x}\left(\frac{\rho_r}{1-\rho_r}\right)^{\sqrt N(g(L)+x)}\prod_{i=0}^{\lfloor LN\rfloor}\left[\sqrt N (g(t_i)+x )\right]_q,
\end{align*} and by $0$ for all other $(x,g,h)\in S_{N}$. Using the expressions for $\rho_{\ell}$ and $\rho_{r}$ (see \eqref{e:rhos}) and \cref{d:scale}, the nonzero part of the weight can also be written as
\begin{align*}
R^{(N)}(x,g,h)=\exp{\left(-2(u+v)x-2vg(L)\right)}\prod_{i=0}^{\lfloor LN\rfloor}\left[\sqrt N (g(t_i)+x )\right]_q.
\end{align*}
We define the corresponding partition function as $\widetilde{\mathcal Z}_{L;u,v}^{(N)}:=N^{-1/2}\sum_{x\in N^{-1/2}\mathbb{Z}_{>0}}\mathbb{E}_{\mathbb{W}_L^{(N)}}\left[R^{(N)} (x,\cdot,\cdot)\right].$
\end{definition}
    In the current notation, which matches \cite{BD23}, it is not true that $R^{(N)}(x,g,h)$ converges to $\Phi_{L;u,v}(x,g,h)$ pointwise. It is first necessary to renormalize the function $R^{(N)}(x,g,h)$. This corresponds to ``recentering'' the weight function to the appropriate height in $x$ in order to obtain the correct pointwise limit. To that end, we define $\Phi^{(N)}_{L;u,v}(x,g,h)$ and $\mathcal{Z}^{(N)}_{L;u,v}$ such that, for appropriate values of $(x,g,h)$, 
    \begin{align}\label{e:varchange}(\mathcal{Z}^{(N)}_{L;u,v})^{-1}\Phi^{(N)}_{L;u,v}(x,g,h)=(\widetilde{\mathcal{Z}}^{(N)}_{L;u,v})^{-1}R^{(N)}(x+\log{(\sqrt{N})},g,h).
\end{align}
\begin{definition} \label{d:hhat} Let $\mathbb{Z}^{(N)}:=\{x\in\mathbb{R}|\sqrt{N}(x+\log{\sqrt{N}})\in\mathbb{Z}_{>0}\}$ and recall that $\mathbb{L}^{(N)}:=N^{-1/2}\sum_{x\in\mathbb{Z}^{(N)}}\delta_{x}$. We define the recentered reference measure on $\mathbb{R}\times C_{0}[0,L]\times C_{0}[0,L]$ by $\mathbb{P}^{(N)}_{L}:=\mathbb{L}^{(N)}\otimes\mathbb{W}_{L}^{(N)}$ (with support $\mathbb{Z}^{(N)}\times W_{[0,L]}^{(N)}$). We define a sequence of functions $\Phi^{(N)}_{L;u,v}:\mathbb{R}\times C_{0}[0,L]\times C_{0}[0,L]\to\mathbb{R}$ by 
    \begin{align*}
    \Phi^{(N)}_{L;u,v}(x,g,h) & := \exp{\left(-2(u+v)x-2vg(L)\right)}\prod_{i=0}^{\lfloor LN\rfloor}\left(1-\frac{e^{-2(g(t_{i})+x)}}{N}\right),
    \end{align*} and note that on the support of $\mathbb{P}_{L}^{(N)}$,
    \begin{align*}
       \Phi^{(N)}_{L;u,v}(x,g,h) = R^{(N)}(x+\log{(\sqrt{N})},g,h)\left(\frac{(1-\rho_{\ell})\rho_{r}}{\rho_{\ell}(1-\rho_{r})}\right)^{-\sqrt{N}\log{(\sqrt{N})}}
    (1-q)^{\lfloor LN\rfloor+1}.
    \end{align*}
    Finally, we define the partition function $\mathcal{Z}^{(N)}_{L;u,v}$ such that $\Phi^{(N)}_{L;u,v}(x,g,h)$ and $\mathcal{Z}^{(N)}_{L;u,v}$ satisfy \eqref{e:varchange}, 
    \begin{align*}\mathcal{Z}^{(N)}_{L;u,v} & :=\mathbb{E}_{\mathbb{P}^{(N)}_{L}}\left[\Phi^{(N)}_{L;u,v} \right].
    \end{align*}
    \end{definition}  
    \cite{BD23} sketches the proof of the pointwise convergence $\Phi^{(N)}_{L;u,v}(x,g,h)\to \Phi_{L;u,v}(x,g,h)$ and conjectures weak convergence of the associated sequence of probability measures. We provide a proof of pointwise convergence in \cref{a:pointwise}.
    \begin{theorem}\label{t:pointwise} Fix $L>0$ and $u,v\in\mathbb{R}$. For all $(x,g,h)\in \mathbb{R}\times C_{0}[0,L]\times C_{0}[0,L]$, as $N\to\infty$ we have $\Phi^{(N)}_{L;u,v}(x,g,h)\to \Phi_{L;u,v}(x,g,h)$ pointwise.\end{theorem}
    \cref{t:pointwise} is purely deterministic. In this article, we demonstrate weak convergence of the associated measures. For the rest of the article, we deal with $\Phi^{(N)}_{L;u,v}(x,g,h)$ instead of $R^{(N)}(x,g,h)$.  
    \begin{definition}\label{d:qqn} We define a sequence of probability measures $\mathbb{Q}^{(N)}_{L;u,v}$ on $\mathbb{R}\times C_{0}[0,L]\times C_{0}[0,L]$ for all $L>0$ and $u,v\in\mathbb{R}$ such that $u+v>0$ through the Radon-Nikodym derivative \begin{align*}
        \frac{d\mathbb{Q}^{(N)}_{L;u,v}}{d\mathbb{P}^{(N)}_{L}}(x,g,h) :=\frac{\Phi^{(N)}_{L;u,v}(x,g,h)}{\mathcal{Z}^{(N)}_{L;u,v}},
    \end{align*} $\mathbb{P}_{L}^{(N)}$-almost everywhere. 
    As a consequence of \cref{t:bldmaintheorem}, the rescaled stationary measure of the open ASEP increment process \eqref{e:rescaledstationary} can be studied through an equality of laws of $$\left\{h^{(N)}_{L;u,v}(x)-h_{L;u,v}^{(N)}(0)\right\}_{x\in[0,L]} \overset{d}{=} \left\{n^{(N)}(x)+m^{(N)}(x)\right\}_{x\in [0,L]}$$ where $(X^{(N)},n^{(N)},m^{(N)})$ are sampled from the measure $\mathbb{Q}^{(N)}_{L;u,v}$. This notation $(X^{(N)},n^{(N)},m^{(N)})$ is not used again. It is intended to suggest the notation in \cref{d:pnotrescaled}. 
\end{definition}
We can now state the main technical result which allows us to prove \cref{t:main1}. 
\begin{theorem}\label{t:main}
    Fix $L>0$ and $u,v\in\mathbb{R}$ such that $u+v>0$. As $N\to\infty$, the measures $\mathbb{Q}^{(N)}_{L;u,v}$ converge weakly to $\mathbb{Q}_{L;u,v}$ (see \cref{d:kpzrn}) on $\mathbb{R}\times C_{0}[0,L]\times C_{0}[0,L]$ where each copy of $C_{0}[0,L]$ is equipped with the uniform topology.
\end{theorem}
The proof of \cref{t:main} is the subject of \cref{s:finite}, \cref{s:conv}, and \cref{s:weakconvergencefinal}. 

\subsection{Proof of \cref{t:main1}}\label{s:measure} 
In this section, we use \cref{t:main} to prove \cref{t:main1}. The version of the open KPZ equation stated in \eqref{e:kpz} takes a different form than that in \cite{CK21,CS18,P19}. We can recover that form of the open KPZ equation (with coefficients of $1/2$ in front of both derivative terms on the right-hand side of the equation and a coefficient of $1$ in front of the noise) by applying the transformation $t\mapsto 2t$ to \eqref{e:kpz}. Since this is only a time transformation, these two versions of the open KPZ equation have the same stationary measure. This justifies applying the results of those papers. The corresponding result is proved in \cite{CS18,P19} on the interval $[0,1]$. We use the following adaptation of \cite[Proposition 3.2]{CK21} to the interval $[0,L]$. The proof carries over after replacing the $N$-site lattice by a $\lfloor LN\rfloor$-site lattice and rescaling the microscopic spatial coordinate accordingly; the constants may depend on the fixed value of $L$.
\begin{proposition}\label{p:old} Consider any $N$-indexed sequence of open ASEPs, which scale under \cref{d:scale}. Assume also that these open ASEPs satisfy
\begin{enumerate}
    \item $4:2:1$ height function scaling: For $t\geq 0$ and $x\in [0,L]$, define
    \begin{align*}
        h^{(N)}_{u,v}(t,x) & := N^{-\frac{1}{2}}h_{N}\left(\frac{e^{\frac{1}{\sqrt{N}}}N^{2}t}{2},x\lfloor LN\rfloor L^{-1}\right) + \frac{Nt}{2}+\frac{t}{24} \\
        z^{(N)}_{u,v}(t,x) & := \exp{\left(h^{(N)}_{u,v}(t,x)\right)}
    \end{align*} where $h_{N}$ is the height function process for the $N$-th open ASEP, whose lattice has $\lfloor LN\rfloor$ sites. 
    \item H\"older bound on initial data: the sequence of initial data $h_{N}(0,\cdot)$ satisfies that for all $\theta\in (0,1/2)$ and $n\in\mathbb{Z}_{>0}$, there exist positive $C(n),C(\theta,n)$ such that for every $x,x'\in [0,L]$ and $N\in\mathbb{Z}_{>0}$, 
    \begin{align*}
        \|z^{(N)}_{u,v}(0,x)\|_{n}\leq C_{L}(n), & & \|z^{(N)}_{u,v}(0,x)-z^{(N)}_{u,v}(0,x')\|_{n}\leq C_{L}(\theta,n)|x-x'|^{\theta},
    \end{align*} where $\|\cdot\|_{n}:=\left(\mathbb{E}[|\cdot|^{n}]\right)^{\frac{1}{n}}$ for expectation taken over $h_{N}(0,\cdot)$.
\end{enumerate}
    Then the sequence of laws of $z^{(N)}_{u,v}\in D([0,T_{0}],C[0,L])$ is tight as $N\to\infty $ for any fixed $T_{0}>0$ and all limit points are in $C([0,T_{0}],C[0,L])$. If there exists a nonnegative $C[0,L]$-valued random variable $z(0,\cdot)$ such that the law of $z^{(N)}_{u,v}(0,\cdot)$ converges weakly to the law of $z(0,\cdot)$ on $C[0,L]$, then the law of $z^{(N)}_{u,v}(\cdot,\cdot)$ converges weakly to that of $z_{u,v}(\cdot,\cdot)$ in $D([0,T_{0}],C[0,L])$ for any $T_{0}>0$, where $z_{u,v}(\cdot,\cdot)$ is the unique mild solution to the stochastic heat equation with boundary parameters compatible with Neumann boundary conditions $u,v\in\mathbb{R}$, and initial data $z(0,\cdot)$. 
\end{proposition}
The H\"{o}lder bounds (2) follow from \cite[Proposition 4.2]{CK21}, whose argument carries over to the $[0,L]$ case. Let $\{N_k\}_{k\geq 1}$ be any subsequence along which the stationary increment profiles converge weakly in $C_{0}[0,L]$ to a random profile. By continuity of the map $f\mapsto e^f$, the corresponding initial Hopf--Cole data also converge, so \cref{p:old} yields convergence along the subsequence $\{N_k\}_{k\geq 1}$ to the stochastic heat equation with corresponding initial data. The finite-$N_{k}$ stationarity identity is preserved under convergence, and using strict positivity we take the logarithm, which allows us to conclude that the law of the increment profile is stationary. Since this holds for any subsequential limit, we have proved the following theorem.
\begin{theorem}\label{t:main2} For any $u,v\in\mathbb{R}$, all subsequential limits of the sequence of open ASEP stationary increment measures under \cref{d:scale} are stationary measures of the open KPZ increment process. 
\end{theorem}
\begin{proof}[Proof of \cref{t:main1}]
By \cref{t:main}, we conclude that $h^{(N)}_{L;u,v}(x)-h^{(N)}_{L;u,v}(0)$ converges weakly to the measure of $\{B_{L}(x)+G_{L;u,v}(x)\}_{x\in [0,L]}$ defined in \cref{t:main1}, where $(X,G_{L;u,v},B_{L})$ are sampled from $\mathbb{Q}_{L;u,v}$. What remains is to conclude that the resulting measure is the stationary measure of the open KPZ increment process on $[0,L]$. This follows from \cref{t:main2}. 
\end{proof}
\section{Two Key Estimates} \label{s:finite}
In this section, we prove two lemmas which we need in \cref{s:conv} and \cref{s:weakconvergencefinal} to establish weak convergence. These are similar to estimates we would need to directly prove tightness of the sequence $\mathbb{Q}^{(N)}_{L;u,v}$. As a consequence of these results, this section also contains a proof that the collection of partition functions $\mathcal{Z}^{(N)}_{L;u,v}$ is uniformly bounded in $N$, and that $\mathcal{Z}_{L;u,v}$ is finite. We begin by stating these lemmas, which are essentially similar results; each provides a bound with a helpful dependence on $x$ for the expectation of the unnormalized weight over the random walk (respectively, Brownian) path measure.

\begin{lemma}\label{l:finitebounds} Fix $L>0$ and $u,v\in\mathbb{R}$ such that $u+v>0$. There exist $N_{0}\in\mathbb{Z}_{>0}$, $\theta>0$ and constants $C_{1},C_{2},C_{3}>0$, which depend on the fixed values $L,u,v$, such that for all $N\geq N_{0}$ and $x\in\mathbb{Z}^{(N)}$,
    \begin{align*}
    \mathbb{E}_{\mathbb{W}^{(N)}_{L}}\left[\Phi^{(N)}_{L;u,v}(x,\cdot,\cdot)\right] &  \leq e^{-2(u+v)x}\left(C_{1}\exp{\left(-\theta W\left(e^{-x}\theta^{-1/2}C_{2}^{1/2}\right)^{2}\right)} + C_{3}N^{-(u+v+1)}\right), 
\end{align*} where $W$ is Lambert's product-log function ($x=W(a)$ solves the equation $xe^x =a$).
\end{lemma}

 \begin{lemma}\label{l:limitbounds}
     Fix $L>0$ and $v\in\mathbb{R}$. There exist $C_{4},C_{5},C_{6},C_{7}>0$ depending on $L,v$ such that for all $x,u\in\mathbb{R}$, if $v\geq 0$, then
\begin{align*}
    \mathbb{E}_{\mathbb{W}_{L}}\left[\Phi_{L;u,v}(x,\cdot,\cdot)\right] \leq e^{-2(u+v)x}C_{4}K_{0}\left(\sqrt{2}e^{-x}\right),
\end{align*} and if $v<0$, then
\begin{align*}
    \mathbb{E}_{\mathbb{W}_{L}}\left[\Phi_{L;u,v}(x,\cdot,\cdot)\right] & \leq e^{-2(u+v)x}C_{5}K_{0}\left(\sqrt{2}e^{-x}\right)  + C_{6}e^{-2(u+v)x+xv}K_{v}\left(C_{7}e^{-x}\right).
\end{align*} In these equations, $K_{v},K_{0}$ denote modified Bessel functions of the
second kind; see \cite{GR}. 
 \end{lemma}
Before proving \cref{l:finitebounds}, we discuss several results which are used in that proof. The first compares a random-walk bridge expectation with the corresponding Brownian-bridge expectation. 
\begin{lemma}\label{l:helpful}
    Fix $L,b>0$. There exist $K:=K(L,b),\lambda:=\lambda(L,b)>0$ such that for every $y:=y(N)\geq 0$ such that $y=O(\log{(N)})$, there exists $N_{0}\in\mathbb{Z}_{>0}$ such that for all $N\geq N_{0}$, $x\in \mathbb{Z}^{(N)}$, and $k\in \llbracket-\lfloor LN\rfloor,\lfloor LN\rfloor\rrbracket$, 
    \begin{align*}
        \mathbb{E}_{\mathbb{W}^{(N)}_{L}}\left[\prod_{i=0}^{\lfloor LN\rfloor}\left(1-\frac{e^{-2(g(t_{i})+x)}}{N}\right)\bigg| g(L)=\frac{k}{\sqrt{N}}\right] &  \leq 4\mathbb{E}_{\mathbb{W}_{L}}\left[\exp{\left(-\frac{e^{-2x}}{2}\int_{0}^{L}e^{-2g(t)}dt\right)}\bigg| g(L)=\frac{k}{\sqrt{N}}\right] 
        \\ & + K\exp{\left(-\lambda y+\frac{bk^2}{\lfloor LN\rfloor}\right)}.
    \end{align*} 
\end{lemma} 
\begin{remark} In the bound above, we ultimately choose $y:=y(N)=2\lambda^{-1}(u+v+1)\log{(\sqrt{N})}$. This choice is important in the proof of \cref{l:finitebounds}.
\end{remark} The other intermediate result that we apply is a translation identity which arises due to translation of Brownian bridge laws. In the proof of \cref{l:finitebounds}, we first use \cref{l:helpful} to pass to the Brownian-bridge setting, because translation by the linear function $mt$ does not preserve the space of discrete random-walk bridge paths.
For $g\in C_0[0,L]$ and $m\in\mathbb{R}$, define the translated path $g_m\in C_0[0,L]$ by $g_m(t):=g(t)-mt$ for $t\in[0,L].$ In particular, $g_m(L)=g(L)-mL.$
\begin{lemma}\label{l:measurechange}
Let $F:C_0[0,L]\to\mathbb{R}$ be a bounded continuous functional. Then, for every $a,m\in\mathbb{R}$,
\begin{align*}
\mathbb{E}_{\mathbb{W}_L}\left[F(g)\bigg|g(L)=a\right]= \mathbb{E}_{\mathbb{W}_L}\left[ F(g_m)\bigg|g(L)=a+mL\right].
\end{align*}
Here $g$ denotes the first coordinate under $\mathbb{W}_L$, and the conditional expectations are understood with respect to the Brownian bridge laws with the indicated endpoints.
\end{lemma}
The proof of this lemma appears in \cref{s:rws}. The following theorem of Dimitrov and Wu~\cite{DW}, adapted to the language of our problem, provides a KMT embedding result for random walk bridges and Brownian bridges, and allows us to prove \cref{l:helpful}.
\begin{theorem}[Bridge KMT coupling; adapted from Theorem 1.2 \cite{DW}]\label{t:KMTbridge} 
    Fix $L>0$. Let $\{X_{\ell}\}_{\ell\geq 1}$ be a sequence of independent, identically distributed random variables satisfying $\mathbb{P}(X_{\ell}=-1)=\mathbb{P}(X_{\ell}=1)=1/4$ and $\mathbb{P}(X_{\ell}=0)=1/2$. For every $b>0$, there exist constants $M,K,\lambda>0$, depending only on $L$ and $b$, such that the following holds. For every $N\in\mathbb{Z}_{>0}$ satisfying $\lfloor LN\rfloor\geq 1$, there exists a probability space (we denote the probability measure on this space by $\mathbb{P}_{\textsc{BKMT}}$, for ``bridge KMT'') carrying, simultaneously for every $k\in\llbracket-\lfloor LN\rfloor,\lfloor LN\rfloor\rrbracket,$ a Brownian bridge  $\widetilde{g}_{k}$ on $[0,L]$ with $\widetilde{g}_{k}(0)=0$ and $\widetilde{g}_{k}(L)=N^{-1/2}k$ and with diffusion coefficient $1/\sqrt{2}$, as well as a continuous process $g_k$ on $[0,L]$. For $ t_i:=iL\lfloor LN\rfloor^{-1},$ for $i\in\llbracket0,\lfloor LN\rfloor\rrbracket,$ the process $g_k$ has the law of the linear interpolation of the process with the law of $N^{-1/2}\sum_{\ell=1}^{i}X_\ell,$ conditioned on $g_{k}(L)=N^{-1/2}k$. 

For every $y\geq0$ and every $k\in\llbracket-\lfloor LN\rfloor,\lfloor LN\rfloor\rrbracket$, the measure $\mathbb{P}_{\textsc{BKMT}}$ satisfies 
\begin{align*}
    \mathbb{P}_{\textsc{BKMT}}\left(\sup_{t\in[0,L]} \left|\widetilde g_k(t)-g_k(t)\right| \geq MN^{-1/2} (\log N+y )\right) 
    \leq K\exp{\left( -\lambda y+\frac{bk^2}{\lfloor LN\rfloor}\right)}.
\end{align*}
The constants $M,K,\lambda$ are independent of $N$, $k$, and $y$.
\end{theorem}
With this result in hand, we can prove \cref{l:helpful}.

\begin{proof}[Proof of \cref{l:helpful}]
    Throughout this proof, we use $\widetilde{g}_{k}$ and $g_{k}$ to denote, respectively, the Brownian bridge path and the random walk bridge path under the BKMT measure. When an expectation is taken with respect to either $\mathbb{W}_{L}^{(N)}$ or $\mathbb{W}_{L}$, we simply use $g$ to denote a path under the appropriate measure. We apply the coupling from  \cref{t:KMTbridge} to write
\begin{align}\label{e:bridgekmt1}
\begin{split}
& \bigg| \mathbb{E}_{\mathbb{W}^{(N)}_{L}}\left[\prod_{i=0}^{\lfloor LN\rfloor }\left(1-\frac{e^{-2(g(t_{i})+x)}}{N}\right)\bigg| g(L)=\frac{k}{\sqrt{N}}\right]-\mathbb{E}_{\mathbb{W}_{L}}\left[\exp{\left(-e^{-2x}\int_{0}^{L}e^{-2g(t)}dt\right)}\bigg| g(L)=\frac{k}{\sqrt{N}}\right] \bigg|
\\ & = \bigg| \mathbb{E}_{\textsc{BKMT}}\left[\prod_{i=0}^{\lfloor LN\rfloor}\left(1-\frac{e^{-2(g_{k}(t_{i})+x)}}{N}\right)-\exp{\left(-e^{-2x}\int_{0}^{L}e^{-2\widetilde{g}_{k}(t)}dt\right)}\right] \bigg|.
\end{split}
\end{align}
We separate the expectation into an event 
\begin{align*}
    A_{k,N}:=\left\{\sup_{t\in[0,L]}\left|\widetilde g_k(t)-g_k(t)\right|<MN^{-\frac{1}{2}}(\log{(N)}+y) \right\},
\end{align*} and its complement.
By \cref{t:KMTbridge}, $\mathbb{P}_{\textsc{BKMT}}\left(A_{k,N}^{c}\right)\leq Ke^{-\lambda y +bk^{2}\lfloor LN\rfloor^{-1}}$. Since both random variables inside the expectation on the right-hand side of \eqref{e:bridgekmt1} are bounded in absolute value by $1$, we conclude that for $K,\lambda>0$ as in \cref{t:KMTbridge}, the expectation over $A_{k,N}^{c}$ is bounded by $2\mathbb{P}_{\textsc{BKMT}}(A_{k,N}^{c})\leq 2Ke^{-\lambda y+bk^{2}\lfloor LN\rfloor^{-1}}$. 
\begin{align*}
     \bigg| \mathbb{E}_{\textsc{BKMT}}\left[\prod_{i=0}^{\lfloor LN\rfloor}\left(1-\frac{e^{-2(g_{k}(t_{i})+x)}}{N}\right)-\exp{\left(-e^{-2x}\int_{0}^{L}e^{-2\widetilde{g}_{k}(t)}dt\right)}\bigg| A^{c}_{k,N}\right] \bigg|\mathbb{P}_{\textsc{BKMT}}\left(A^{c}_{k,N}\right) \leq 2Ke^{-\lambda y+bk^{2}\lfloor LN\rfloor^{-1}}.
 \end{align*}
It remains to bound the conditional expectation over $A_{k,N}$,
\begin{align}\label{e:KMTAterms}\begin{split}
    & \bigg| \mathbb{E}_{\textsc{BKMT}}\left[\prod_{i=0}^{\lfloor LN\rfloor}\left(1-\frac{e^{-2(g_{k}(t_{i})+x)}}{N}\right)-\exp{\left(-e^{-2x}\int_{0}^{L}e^{-2\widetilde{g}_{k}(t)}dt\right)}\bigg| A_{k,N} \right] \bigg|\mathbb{P}_{\textsc{BKMT}}(A_{k,N})
    \\ & \leq \bigg| \mathbb{E}_{\textsc{BKMT}}\left[\prod_{i=0}^{\lfloor LN\rfloor}\left(1-\frac{e^{-2(g_{k}(t_{i})+x)}}{N}\right)-\exp{\left(-e^{-2x}\int_{0}^{L}e^{-2g_{k}(t)}dt\right)}\bigg| A_{k,N}\right] \bigg|\mathbb{P}_{\textsc{BKMT}}(A_{k,N})
    \\ & + \bigg| \mathbb{E}_{\textsc{BKMT}}\left[\exp{\left(-e^{-2x}\int_{0}^{L}e^{-2g_{k}(t)}dt\right)}-\exp{\left(-e^{-2x}\int_{0}^{L}e^{-2\widetilde{g}_{k}(t)}dt\right)}\bigg| A_{k,N} \right] \bigg|\mathbb{P}_{\textsc{BKMT}}(A_{k,N}).
    \end{split}
\end{align}
 We define the notation $\delta:=\delta(N)= MN^{-1/2}\left(\log{(N)}+y\right)$. On the event $A_{k,N}$, $e^{-2\delta}e^{-2\widetilde g_k(t)}\leq e^{-2g_k(t)}\leq e^{2\delta}e^{-2\widetilde g_k(t)}.$ Consequently,
\begin{align*}
    e^{-2\delta}e^{-2x}\int_0^L e^{-2\widetilde g_k(t)}dt\leq e^{-2x}\int_0^L e^{-2g_k(t)}dt\leq e^{2\delta}e^{-2x}\int_0^L e^{-2\widetilde g_k(t)}dt.
\end{align*}
Therefore, the second term on the right-hand side of \eqref{e:KMTAterms} is bounded by 
\begin{align*}
    & \leq \mathbb{E}_{\textsc{BKMT}}\left[\exp{\left(-e^{-2x}\int_{0}^{L}e^{-2\widetilde{g}_{k}(t)}dt\right)}\max{\left\{\begin{matrix}1-
\exp{\left(
-\left(e^{2\delta}-1\right)e^{-2x}
\int_0^L e^{-2\widetilde g_k(t)}dt
\right)},\\ \exp{\left(
\left(1-e^{-2\delta}\right)e^{-2x}
\int_0^L e^{-2\widetilde g_k(t)}dt
\right)}-1\end{matrix}\right\}} \mathbf{1}_{A_{k,N}} \right].
\end{align*} We consider the $(0,\infty)$-valued random variable $X_{k}=e^{-2x}\int_{0}^{L}e^{-2\widetilde{g}_{k}(t)}dt$ (we suppress the dependence on $N,x$ in this notation). By the joint-density formula of Matsumoto and Yor \cite[Theorem 4.1]{MY05}, after the appropriate sign and scaling changes, under $\mathbb{W}_L$ the pair
\begin{align*} \left(e^{-2x}\int_0^L e^{-2g(t)}dt,g(L)\right)
\end{align*} has a joint density on $(0,\infty)\times\mathbb{R}$. Since $g(L)$ has a Gaussian density, the conditional law of the first coordinate given $g(L)=kN^{-1/2}$ has an explicit density. Since the process $\widetilde g_k$ has the law of the first coordinate under $\mathbb{W}_{L}$ conditioned on $g(L)=kN^{-1/2}$ under $\mathbb{P}_{\textsc{BKMT}}$, it follows that $X_k$ has the corresponding conditional density, which we denote by $\rho_{X,k}$ (we continue to suppress the dependence on $N,x$). From this, we deduce that the preceding expression is bounded above by
\begin{align*}
    & \leq \int_{0}^{\infty}e^{-r}\max{\left\{1-e^{-r(e^{2\delta}-1)},e^{-r(e^{-2\delta}-1)}-1\right\}}\rho_{X,k}(r)dr.
\end{align*}
Since $y=O(\log N)$, as $N\to\infty$, $\delta=\delta(N)=MN^{-1/2} (\log N+y )\to 0$. Therefore, for $N$ sufficiently large, i.e. after increasing $N_0$ if necessary, for every $N\geq N_0$, $0<e^{2\delta}-1\leq 1/4$ and $0<1-e^{-2\delta}\leq 1/4$. This implies that the expression above is bounded by
\begin{align*}
     & \leq \int_{0}^{\infty}e^{-r}\max{\left\{ 1-e^{-r/4} ,e^{r/4}-1 \right\}}\rho_{X,k}(r)dr.
\end{align*}
For every $r\geq 0$, $e^{-r}(1-e^{-r/4}),e^{-r}(e^{r/4}-1)\leq e^{-r/2}.$ Consequently, $\max{\{1-e^{-r/4},e^{r/4}-1\}}\leq e^{r/2}.$ Therefore, the preceding integral is bounded by
\begin{align*}
\int_0^\infty e^{-r/2}\rho_{X,k}(r)dr 
& =\mathbb{E}_{\mathbb{W}_L}\left[\exp{\left(-\frac{e^{-2x}}{2}\int_0^L e^{-2g(t)}dt\right)}\bigg|g(L)=\frac{k}{\sqrt N}\right].
\end{align*}This completes the bound on the second term on the right-hand side
of \eqref{e:KMTAterms}.
We now bound the first term on the right-hand side of \eqref{e:KMTAterms}. The definition of the linearly interpolated random-walk path implies
\begin{align*}
\int_0^L e^{-2(g_k(t)+x)} dt
& = \frac{L}{\lfloor LN\rfloor} \sum_{i=0}^{\lfloor LN\rfloor-1}e^{-2(g_k(t_i)+x)} \int_0^1 \exp{\left( -2s(g_k(t_{i+1})-g_k(t_i)) \right)} ds.
\end{align*}
Since each increment of the first coordinate of the rescaled random-walk path has absolute value at most $N^{-1/2}$, it follows that $|g_k(t_{i+1})-g_k(t_i)|\leq N^{-1/2}.$ Consequently, for every $s\in[0,1]$, $e^{-2s(g_k(t_{i+1})-g_k(t_i))}\leq e^{2/\sqrt N},$ and, therefore,
\begin{align*}
e^{-2/\sqrt N}\frac{\lfloor LN\rfloor}{LN}\int_0^L e^{-2(g_k(t)+x)}dt \leq  \frac{1}{N}\sum_{i=0}^{\lfloor LN\rfloor}e^{-2(g_k(t_i)+x)}.
\end{align*}
For $x\in\mathbb{Z}^{(N)}$, the discrete product either vanishes or
every factor in the product belongs to $(0,1)$. Therefore, since $1-z\leq e^{-z}$ when the product is nonzero, we conclude that
\begin{align*}
\prod_{i=0}^{\lfloor LN\rfloor} \left(1-\frac{e^{-2(g_k(t_i)+x)}}{N}\right) 
\leq\exp{\left(-\frac{\lfloor LN\rfloor}{LN}e^{-2/\sqrt N}\int_0^L e^{-2(g_k(t)+x)}dt\right)},
\end{align*}
Since $0<\lfloor LN\rfloor(LN)^{-1}e^{-2/\sqrt N}\leq 1,$
It follows that
\begin{align*}
\exp{\left(-e^{-2x}\int_0^L e^{-2g_k(t)}dt\right)}
& \leq\exp{\left(-\frac{\lfloor LN\rfloor}{LN}e^{-2/\sqrt N}\int_0^L e^{-2(g_k(t)+x)}dt\right)}.
\end{align*}
Therefore, the absolute difference of the two terms appearing in the first term in \eqref{e:KMTAterms} is also bounded by the same quantity. Therefore, on the event $A_{k,N}$,
\begin{align*}
e^{-2\delta}\int_0^L e^{-2(\widetilde g_k(t)+x)}dt\leq \int_0^L e^{-2(g_k(t)+x)}dt.
\end{align*}
Moreover, $\lfloor LN\rfloor(LN)^{-1}e^{-2/\sqrt N}e^{-2\delta}\to 1$ as $N\to\infty$. Therefore, after possibly increasing $N_0$, we conclude that for every $N\geq N_0$, $\lfloor LN\rfloor(LN)^{-1}e^{-2/\sqrt N}e^{-2\delta}\geq 1/2$. Therefore, on the event $A_{k,N}$,
\begin{align*}
\exp{\left(-\frac{\lfloor LN\rfloor}{LN}e^{-2/\sqrt N}\int_0^L e^{-2(g_k(t)+x)}dt\right)} \leq\exp{\left(-\frac{e^{-2x}}2\int_0^L e^{-2\widetilde g_k(t)}dt\right)}.
\end{align*} Finally, we conclude that the first term on the right-hand side of \eqref{e:KMTAterms} is bounded by
\begin{align*}
\mathbb{E}_{\textsc{BKMT}}\left[\exp{\left(-\frac{e^{-2x}}2\int_0^L e^{-2\widetilde g_k(t)}dt\right)}\mathbf{1}_{A_{k,N}}\right] 
& \leq \mathbb{E}_{\mathbb{W}_L}\left[\exp{\left(-\frac{e^{-2x}}2\int_0^L e^{-2g(t)}dt\right)}\bigg|g(L)=\frac{k}{\sqrt N}\right].
\end{align*} 
Combining the bounds on the two terms in \eqref{e:KMTAterms} with the contribution from $A_{k,N}^{c}$, we conclude that the absolute difference in \eqref{e:bridgekmt1} is bounded by
\begin{align*}
& 2\mathbb{E}_{\mathbb{W}_L}\left[\exp{\left( -\frac{e^{-2x}}2 \int_0^L e^{-2g(t)}dt \right)} \bigg| g(L)=\frac{k}{\sqrt N} \right] + 2K\exp{\left(-\lambda y+\frac{bk^2}{\lfloor LN\rfloor}\right)}.
\end{align*}
Moreover,
\begin{align*}
\mathbb{E}_{\mathbb{W}_L} \left[\exp{\left(-e^{-2x}\int_0^L e^{-2g(t)}dt\right)} \bigg| g(L)=\frac{k}{\sqrt N} \right] 
\leq \mathbb{E}_{\mathbb{W}_L}\left[\exp{\left(-\frac{e^{-2x}}2\int_0^L e^{-2g(t)} dt\right)} \bigg| g(L)=\frac{k}{\sqrt N} \right].
\end{align*}
Therefore, by combining our bounds and enlarging the constant $K$, we obtain the result.
\end{proof}
We will also need a comparison in the reverse direction. We claim that, uniformly for $m$ in any fixed compact subset of $\mathbb{R}$, after possibly increasing $N_0$,
\begin{align}\label{e:reverse}
\begin{split}
&
\mathbb{E}_{\mathbb{W}_L}\left[\exp{\left(-\frac{e^{-2x}}{2}\int_0^L e^{-2(g(t)-mt)} dt\right)} \bigg|g(L)=\frac{k}{\sqrt N}\right]
\\ & \leq 4\mathbb{E}_{\mathbb{W}_L^{(N)}}\left[\prod_{i=0}^{\lfloor LN\rfloor}\left(1+\frac{e^{-2(g(t_i)-mt_i+x)}}{8N}\right)^{-1}\bigg|g(L)=\frac{k}{\sqrt N}\right]+K\exp{\left(-\lambda y+\frac{bk^2}{\lfloor LN\rfloor}\right)}.
\end{split}
\end{align}
The proof uses the same bridge-KMT coupling and the same decomposition into a good coupling event and its complement as the proof of \cref{l:helpful}.

We will also use  a binomial-shift estimate in the proof of \cref{l:finitebounds}. Its proof is given in \cref{s:rws}.
\begin{proposition}\label{l:combapprox}
Fix $L>0$ and $a\in\mathbb{R}$. There exist $D_{L,a}>0$ and $N_{L,a}\in\mathbb{Z}_{>0}$ which may depend on $a,L$ such that, for every $N\geq N_{L,a}$ and every $k\in\mathbb{Z}$ satisfying $|k|<\lfloor LN\rfloor^{5/6},$ and $|k+\lfloor a\lfloor LN\rfloor N^{-1/2}\rfloor|<\lfloor LN\rfloor^{5/6},$ we have
\begin{align*}
\exp{\left(-\frac{2ak}{\sqrt N}\right)}\binom{2\lfloor LN\rfloor}{\lfloor LN\rfloor+k}\leq D_{L,a}\binom{2\lfloor LN\rfloor}{\lfloor LN\rfloor+k+\lfloor a\lfloor LN\rfloor N^{-1/2}\rfloor 
}.
\end{align*}
\end{proposition}
\begin{remark}
    The choice of $N^{c}$ with $c=5/6$ as the threshold in \cref{l:combapprox} is arbitrary: all that is necessary for the lemma to hold is that $c\leq 5/6$, and all that is necessary for this lemma to be used in the ongoing proof is that $c>3/4$.
\end{remark}

The preceding comparisons in \cref{l:helpful} and \eqref{e:reverse} require only that $MN^{-1/2} (\log N+y )\to 0.$ In the proof of \cref{l:finitebounds}, we choose $y$ as in the preceding remark to control the KMT error. These comparisons allow us to pass to Brownian bridges, apply the translation identity in \cref{l:measurechange}, and return to a positive discrete functional whose conditional expectations can be recombined over the endpoint distribution. We can now begin the proof of \cref{l:finitebounds}.
\begin{proof}[Proof of \cref{l:finitebounds}]
Conditioning on $g(L)$, we find that
\begin{multline}\label{e:conditioned}
\mathbb{E}_{\mathbb{W}^{(N)}_{L}}\left[\Phi^{(N)}_{L;u,v}(x,g,h)\right] = e^{-2(u+v)x}4^{-\lfloor LN\rfloor}
\\   \cdot \sum_{k=-\lfloor LN\rfloor}^{\lfloor LN\rfloor}\exp{\left(-\frac{2vk}{\sqrt N}\right)}\binom{2\lfloor LN\rfloor}{\lfloor LN\rfloor+k}  \mathbb{E}_{\mathbb{W}^{(N)}_{L}}\left[\prod_{i=0}^{\lfloor LN\rfloor}\left(1-\frac{e^{-2(g(t_{i})+x)}}{N}\right) \bigg|g(L)=\frac{k}{\sqrt N}\right].
\end{multline}
Fix $b>0$, to be chosen sufficiently small below. Applying \cref{l:helpful} to \eqref{e:conditioned}, and then applying \cref{l:measurechange} with translation parameter $m := \lfloor v\lfloor LN\rfloor N^{-1/2}\rfloor L^{-1}N^{-1/2},$ we obtain
\begin{align*}
& \mathbb{E}_{\mathbb{W}^{(N)}_{L}}\left[\Phi^{(N)}_{L;u,v}(x,g,h)\right] \leq e^{-2(u+v)x}4^{-\lfloor LN\rfloor+1}\sum_{k=-\lfloor LN\rfloor}^{\lfloor LN\rfloor}e^{-2vk/\sqrt N}\binom{2\lfloor LN\rfloor}{\lfloor LN\rfloor+k}
\\ & \cdot\mathbb{E}_{\mathbb{W}_{L}}\left[\exp{\left(-\frac{e^{-2x}}{2}\int_0^L e^{-2\left(g(t)-mt\right)}dt\right)}\bigg|g(L)=\frac{k}{\sqrt N} +mL\right]
\\ &+K e^{-2(u+v)x-\lambda y}4^{-\lfloor LN\rfloor}\sum_{k=-\lfloor LN\rfloor}^{\lfloor LN\rfloor}\exp{\left(-\frac{2vk}{\sqrt N}+\frac{bk^2}{\lfloor LN\rfloor}\right)}\binom{2\lfloor LN\rfloor}{\lfloor LN\rfloor+k}.
\end{align*}
Before applying \eqref{e:reverse}, we split the first sum into the set of $k$ such that $|k|<\lfloor LN\rfloor^{5/6}$ and $|k+mL\sqrt{N}|<\lfloor LN\rfloor^{5/6}$ are both satisfied, and the complement of this set. On the complement, the conditional expectation is bounded by $1$. It follows that $m\to v$ as $N\to\infty$. Thus, for all sufficiently large $N$, $m$ lies in a fixed compact subset of $\mathbb{R}$, and the estimate of \eqref{e:reverse} applies uniformly.

Thus, for all $k$ satisfying $|k|<\lfloor LN\rfloor^{5/6}$ and $|k+mL\sqrt{N}|<\lfloor LN\rfloor^{5/6}$, we apply \eqref{e:reverse} with the endpoint index $k+mL\sqrt{N}$ to obtain
\begin{align*}
&
\mathbb{E}_{\mathbb{W}_L}\left[\exp{\left(-\frac{e^{-2x}}2\int_0^L e^{-2(g(t)-mt)} dt\right)}\bigg|g(L)=\frac{k}{\sqrt N} +mL\right]
\\ & \leq 4\mathbb{E}_{\mathbb{W}_L^{(N)}}\left[\prod_{i=0}^{\lfloor LN\rfloor}\left(1+\frac{e^{-2(g(t_i)-mt_i+x)}}{8N}\right)^{-1} \bigg|g(L)=\frac{k}{\sqrt N}+mL\right] +K\exp{\left(-\lambda y+\frac{b\left(k+mL\sqrt{N}\right)^2}{\lfloor LN\rfloor}\right)}.
\end{align*}
Substituting these bounds into the first sum above, and increasing $K$ if necessary, we obtain
\begin{align}\label{e:newconditions}\begin{split}
\mathbb{E}_{\mathbb{W}_L^{(N)}} \left[\Phi_{L;u,v}^{(N)}(x,g,h)\right]& \leq\frac{e^{-2(u+v)x}}{4^{\lfloor LN\rfloor-2}}\sum_{\substack{-\lfloor LN\rfloor\leq k\leq\lfloor LN\rfloor \\ |k|<\lfloor LN\rfloor^{5/6} \\ \left|k+mL\sqrt{N}\right|<\lfloor LN\rfloor^{5/6}}}\exp{\left(-\frac{2vk}{\sqrt N}\right)}\binom{2\lfloor LN\rfloor}{\lfloor LN\rfloor+k}
\\ & \cdot \mathbb{E}_{\mathbb{W}_L^{(N)}}\left[\prod_{i=0}^{\lfloor LN\rfloor}\left(1+\frac{e^{-2(g(t_i)-mt_i+x)}}{8N}\right)^{-1}\bigg|g(L)=\frac{k}{\sqrt N}+mL\right]
\\ & + e^{-2(u+v)x}\left(\textsc{ERR}_{1}+\textsc{ERR}_{2}+\textsc{ERR}_{3}\right).
\end{split}
\end{align}
Where we have defined 
\begin{align*} 
\textsc{ERR}_{1} & :=\frac{K e^{-\lambda y}}{4^{\lfloor LN\rfloor-1}}\sum_{\substack{-\lfloor LN\rfloor\leq k\leq\lfloor LN\rfloor\\ |k|<\lfloor LN\rfloor^{5/6}\\ |k+mL\sqrt{N}|<\lfloor LN\rfloor^{5/6}}}\exp{\left(-\frac{2vk}{\sqrt N}+\frac{b\left(k+mL\sqrt{N}\right)^2}{\lfloor LN\rfloor}\right)} \binom{2\lfloor LN\rfloor}{\lfloor LN\rfloor+k},
\\ \textsc{ERR}_{2} & :=\frac{1}{4^{\lfloor LN\rfloor-1}}\sum_{\substack{-\lfloor LN\rfloor\leq k\leq\lfloor LN\rfloor\\ |k|\geq\lfloor LN\rfloor^{5/6} \mathrm{ or }\\ |k+mL\sqrt{N}|\geq\lfloor LN\rfloor^{5/6}}}\exp{\left(-\frac{2vk}{\sqrt N}\right)}\binom{2\lfloor LN\rfloor}{\lfloor LN\rfloor+k},
\\ \textsc{ERR}_{3} & :=\frac{K e^{-\lambda y}}{4^{\lfloor LN\rfloor}}\sum_{k=-\lfloor LN\rfloor}^{\lfloor LN\rfloor}\exp{\left(-\frac{2vk}{\sqrt N}+\frac{bk^2}{\lfloor LN\rfloor}\right)}\binom{2\lfloor LN\rfloor}{\lfloor LN\rfloor+k}.
\end{align*}

We begin by estimating the first term of \eqref{e:newconditions}. Applying \cref{l:combapprox} with $a=v$, for all $k$ which appear in that sum we have
\begin{align*}
\exp{\left(-\frac{2vk}{\sqrt N}\right)}\binom{2\lfloor LN\rfloor}{\lfloor LN\rfloor+k} 
\leq D_{L,v}\binom{2\lfloor LN\rfloor}{\lfloor LN\rfloor+k+mL\sqrt{N}}.
\end{align*}
We make the change of variables $j=k+mL\sqrt{N}$ and use the fact that the product which appears inside the conditional expectation is always positive to enlarge the resulting sum to every $j\in\llbracket-\lfloor LN\rfloor,\lfloor LN\rfloor\rrbracket.$ Finally we can sum over $j$ to find that the first term in \eqref{e:newconditions} is
bounded by
\begin{align*}
& \frac{D_{L,v}e^{-2(u+v)x}}{4^{\lfloor LN\rfloor-2}}\sum_{j=-\lfloor LN\rfloor}^{\lfloor LN\rfloor}\binom{2\lfloor LN\rfloor}{\lfloor LN\rfloor+j}\mathbb{E}_{\mathbb{W}_L^{(N)}}\left[\prod_{i=0}^{\lfloor LN\rfloor}\left(1+\frac{e^{-2(g(t_i)-mt_i+x)}}{8N}\right)^{-1}\bigg|g(L)=\frac{j}{\sqrt N}\right] 
\\ & =16D_{L,v}e^{-2(u+v)x}\mathbb{E}_{\mathbb{W}_L^{(N)}}\left[\prod_{i=0}^{\lfloor LN\rfloor}\left(1+\frac{e^{-2(g(t_i)-mt_i+x)}}{8N}\right)^{-1}\right].
\end{align*}
We proceed to estimate the subsequent terms. We fix $b>0$ sufficiently small; for example, take $0<b\leq 1/16$. Due to Stirling's formula, there exists a constant $C>0$ such that for every $k\in\llbracket-\lfloor LN\rfloor,\lfloor LN\rfloor\rrbracket$,
\begin{align*}
\frac{1}{4^{\lfloor LN\rfloor}}\binom{2\lfloor LN\rfloor}{\lfloor LN\rfloor+k}
\leq\frac{C}{\sqrt{\lfloor LN\rfloor}}\exp{\left(-\frac{k^2}{2\lfloor LN\rfloor}\right)}
\end{align*}
Furthermore, since $2|v||k|N^{-1/2}\leq k^2(8\lfloor LN\rfloor)^{-1}+8Lv^2$ and $|mL\sqrt{N}|\leq |v|\lfloor LN\rfloor N^{-1/2}+1,$ there is also a constant $C_{L,v}>0$ such that
\begin{align*}
b\left(k+mL\sqrt{N}\right)^2\leq k^{2}/8+ C_{L,v}\lfloor LN\rfloor.
\end{align*}
It follows that
\begin{align*}
\frac{1}{4^{\lfloor LN\rfloor}}\exp{\left(-\frac{2vk}{\sqrt N}+\frac{b\left(k+mL\sqrt{N}\right)^2}{\lfloor LN\rfloor}\right)}\binom{2\lfloor LN\rfloor}{\lfloor LN\rfloor+k}\leq\frac{C_{L,v}}{\sqrt{\lfloor LN\rfloor}}\exp{\left(-\frac{k^2}{4\lfloor LN\rfloor}\right)},
\end{align*}
and similarly
\begin{align*}
& \frac{1}{4^{\lfloor LN\rfloor}}\exp{\left(-\frac{2vk}{\sqrt N}+\frac{bk^2}{\lfloor LN\rfloor}\right)}\binom{2\lfloor LN\rfloor}{\lfloor LN\rfloor+k}
\leq \frac{C_{L,v}}{\sqrt{\lfloor LN\rfloor}}\exp{\left(-\frac{k^2}{4\lfloor LN\rfloor}\right)}.
\end{align*}
Therefore, since  the sum 
\begin{align*}
\frac{1}{\sqrt{\lfloor LN\rfloor}}\sum_{k=-\lfloor LN\rfloor}^{\lfloor LN\rfloor} \exp{\left(-\frac{k^2}{4\lfloor LN\rfloor}\right)}
\end{align*}
is uniformly bounded in $N$, we conclude that after absorbing the
factors $K$ and $4$ into the constant, the expressions $\textsc{ERR}_{1},\textsc{ERR}_{3}$ are both bounded as
\begin{align*}
    \textsc{ERR}_{1},\textsc{ERR}_{3} \leq C_{L,v}e^{-\lambda y}.
\end{align*}
We now turn to producing a bound on $\textsc{ERR}_{2}$. For all sufficiently large $N$, $|mL\sqrt{N}|\leq \lfloor LN\rfloor^{5/6}/2.$
Consequently, if either $|k|\geq \lfloor LN\rfloor^{5/6}$ or $|k+mL\sqrt{N}|\geq \lfloor LN\rfloor^{5/6},$ then $|k|\geq \lfloor LN\rfloor^{5/6}/2$. Thus, applying previous estimates, we can conclude that there exists $c>0$ such that
\begin{align*}
\textsc{ERR}_{2}\leq\frac{C_{L,v}}{\sqrt{\lfloor LN\rfloor}}\sum_{\substack{|k|\geq\lfloor LN\rfloor^{5/6}/2 \\|k|\leq\lfloor LN\rfloor}}\exp{\left(-\frac{k^2}{4\lfloor LN\rfloor}\right)} \leq C_{L,v}\exp{\left(-c\lfloor LN\rfloor^{2/3}\right)}.
\end{align*}

We define the notation  $\textsc{ERR}_{4}:=\textsc{ERR}_{1}+\textsc{ERR}_{2}+\textsc{ERR}_{3}.$ Now, when we choose $y=2\lambda^{-1}(u+v+1)\log(\sqrt N),$ then $e^{-\lambda y}=N^{-(u+v+1)}.$ Furthermore, for all $N$ sufficiently large,  $\exp{\left( -c\lfloor LN\rfloor^{2/3} \right)} \leq e^{-\lambda y}.$ It follows that $\textsc{ERR}_{4} \leq C_{L,u,v}e^{-\lambda y}.$

Thus, after relabeling constants, there exist $D_1,D_2>0$, depending only on $L,u,v$, such that \eqref{e:newconditions} is bounded by
\begin{align}\label{e:resummedbound}
\mathbb{E}_{\mathbb{W}_L^{(N)}}\left[\Phi_{L;u,v}^{(N)}(x,g,h)\right] &\leq D_1 e^{-2(u+v)x}\mathbb{E}_{\mathbb{W}_L^{(N)}}\left[\prod_{i=0}^{\lfloor LN\rfloor}\left(1+\frac{e^{-2(g(t_i)-mt_i+x)}}{8N}\right)^{-1}\right] + D_2e^{-2(u+v)x-\lambda y}.
\end{align}

To continue, for $s>0$, we define 
\begin{align*}
    A_{N,s}:=\left\{\max_{0\leq i\leq\lfloor LN\rfloor}|g(t_i)|>s\right\}.
\end{align*}
Using this event, we make the decomposition
\begin{align*}
\mathbb{E}_{\mathbb{W}_L^{(N)}}\left[\prod_{i=0}^{\lfloor LN\rfloor}\left(1+\frac{e^{-2(g(t_i)-mt_i+x)}}{8N}\right)^{-1}\right]
& =\mathbb{E}_{\mathbb{W}_L^{(N)}}\left[\prod_{i=0}^{\lfloor LN\rfloor}\left(1+\frac{e^{-2(g(t_i)-mt_i+x)}}{8N}\right)^{-1}\bigg|A_{N,s}\right]\mathbb{W}_L^{(N)}(A_{N,s})
\\ & +\mathbb{E}_{\mathbb{W}_L^{(N)}}\left[\prod_{i=0}^{\lfloor LN\rfloor}\left(1+\frac{e^{-2(g(t_i)-mt_i+x)}}{8N}\right)^{-1}\bigg|A_{N,s}^c\right]\mathbb{W}_L^{(N)}(A_{N,s}^c).
\end{align*}
Since the product inside the conditional expectation is bounded by $1$, the first term is bounded by $\mathbb{W}_L^{(N)}(A_{N,s}).$
By \cite[Proposition 2.1.2]{L10}, there exist constants $C,\theta>0$, depending only on $L$, such that, for every $s>0$, $\mathbb{W}_{L}^{(N)}(A_{N,s})\leq Ce^{-\theta s^2}.$

For all sufficiently large $N$, the value of $m$ chosen above satisfies $|m|\leq |v|+1.$ On the event $A_{N,s}^c$, for every $i\in\llbracket0,\lfloor LN\rfloor\rrbracket$, we therefore have $g(t_i)-mt_i \leq s+(|v|+1)L.$ Consequently,
\begin{align*} e^{-2(g(t_i)-mt_i+x)} \geq e^{-2(x+s+(|v|+1)L)},
\end{align*}
and, therefore, on $A_{N,s}^c$,
\begin{align*}
\prod_{i=0}^{\lfloor LN\rfloor}\left(1+\frac{e^{-2(g(t_i)-mt_i+x)}}{8N}\right)^{-1} & \leq \left(1+\frac{e^{-2(x+s+(|v|+1)L)}}{8N}\right)^{-\lfloor LN\rfloor-1}.
\end{align*}
Since $x\in\mathbb{Z}^{(N)}$, we have $e^{-2x}/N<1$, and therefore the expression inside the parentheses on the right hand side is bounded above by $9/8$. Using the fact that  $z/2\leq \log(1+z)$ for $z\in [0,1]$, we therefore conclude
\begin{align*}
\prod_{i=0}^{\lfloor LN\rfloor}\left(1+\frac{e^{-2(g(t_i)-mt_i+x)}}{8N}\right)^{-1} & \leq\exp{\left(-\frac{\lfloor LN\rfloor+1}{16N}e^{-2(x+s+(|v|+1)L)}\right)} \leq\exp{\left(-C_2e^{-2(x+s)}\right)},
\end{align*}
where $C_2 :=Le^{-2(|v|+1)L}/16$. Therefore,
\begin{align}\label{e:psimaxbound}
\mathbb{E}_{\mathbb{W}_L^{(N)}}\left[\prod_{i=0}^{\lfloor LN\rfloor}\left(1+\frac{e^{-2(g(t_i)-mt_i+x)}}{8N}\right)^{-1}\right]\leq Ce^{-\theta s^2}+\exp{\left(-C_2e^{-2(x+s)}\right)}.
\end{align}

We balance the two exponential terms by choosing $s>0$ such that $\theta s^2=C_2e^{-2(x+s)}$. Solving for $s$ yields $s=W\left(e^{-x}\theta^{-1/2}C_2^{1/2}\right).$ Therefore, 
\begin{align*}
\mathbb{E}_{\mathbb{W}_L^{(N)}}\left[\prod_{i=0}^{\lfloor LN\rfloor}\left(1+\frac{e^{-2(g(t_i)-mt_i+x)}}{8N}\right)^{-1}\right] & \leq (C+1)\exp{\left(-\theta W\left(e^{-x}\theta^{-1/2}C_2^{1/2}\right)^2\right)}.
\end{align*}

Substituting this estimate and our chosen value of $y$ into \eqref{e:resummedbound} results in 
\begin{align*}
\mathbb{E}_{\mathbb{W}_L^{(N)}}\left[\Phi_{L;u,v}^{(N)}(x,g,h)\right] & \leq D_1 (C+1)e^{-2(u+v)x} \exp{\left(-\theta W\left(e^{-x}\theta^{-1/2}C_2^{1/2}\right)^2\right)}+D_2 e^{-2(u+v)x}N^{-(u+v+1)}.
\end{align*}
Setting $C_1:=D_1(C+1)$ and $C_3:=D_2$ finishes the proof.
\end{proof}
As a consequence of \cref{l:finitebounds}, the sequence $\{\mathcal{Z}^{(N)}_{L;u,v}\}_{N\in\mathbb{Z}_{>0}}$ is uniformly bounded. 
\begin{corollary}\label{c:part}
Fix $L>0$ and $u,v\in\mathbb{R}$ such that $u+v>0$. There exists $C>0$, depending only on $L,u,v$, such that for all $N\in\mathbb{Z}_{>0}$ such that $\lfloor LN\rfloor \geq 1$, $\mathcal Z_{L;u,v}^{(N)}\leq C$. 
\end{corollary}

This corollary follows by writing $\mathcal Z_{L;u,v}^{(N)}$ as a sum of terms of the form which appear in \cref{l:finitebounds} over $x\in\mathbb{Z}^{(N)}$, and then applying \cref{l:finitebounds}. The contribution of the second term in that bound becomes $N^{-3/2}\sum_{j=1}^{\infty} e^{-2(u+v)j/\sqrt N},$ which is uniformly bounded in $N$. For the first term, when $x\geq0$, the summand is bounded by $e^{-2(u+v)x}$, while, as $x\to-\infty$, $W\left(e^{-x}\theta^{-1/2}C_2^{1/2}\right)\sim |x|,$ and so the summand is bounded by $e^{-c x^2}$ for some $c>0$. The preceding argument gives a uniform bound for $N\geq N_0$. For each of the finitely many $N<N_0$ satisfying $\lfloor LN\rfloor\geq1$, the discrete product belongs to $[0,1]$, and thus
\begin{align*}
\mathcal Z_{L;u,v}^{(N)}\leq N^{-1/2}\mathbb{E}_{\mathbb{W}_L^{(N)}}\left[e^{-2vg(L)}\right]\sum_{x\in\mathbb{Z}^{(N)}}e^{-2(u+v)x}<\infty.
\end{align*}
The sum is finite because $u+v>0$. Increasing the constant allows us to obtain the claimed uniform bound for every $N<N_{0}$ satisfying $\lfloor LN\rfloor\geq1$.

Having now finished the proof of \cref{l:finitebounds}, we proceed to prove \cref{l:limitbounds}

\begin{proof}[Proof of \cref{l:limitbounds}]
Let $g$ denote the first Brownian coordinate under $\mathbb{W}_L$. By the corollary in Section 6.2 of \cite{Y92}; see also \cite[Theorem 4.1 and (2.9)]{MY05}, for nonnegative Borel measurable functions $f,\varphi:(0,\infty)\to[0,\infty)$, we have
\begin{align}\label{e:YOR}\begin{split}
 \mathbb{E}_{\mathbb{W}_L}\left[f\left(e^{-g(L)}\right)\varphi\left(\int_0^L e^{-2g(t)} dt\right)\right]  &  = D_1 \int_0^\infty\int_0^\infty f(y) \varphi\left(\frac{2}{r}\right)\exp{\left(-\frac{r}{2}(1+y^2)\right)} 
 \\ & \cdot\left(\int_0^\infty \exp{\left(-\frac{p^2}{L}-yr\cosh(p)\right)}\sinh(p)\sin{\left(\frac{2\pi p}{L}\right)}dp\right)dr dy,
\end{split}
\end{align}
for a constant $D_1=D_1(L)>0$. After multiplying by $e^{-2(u+v)x}$, we obtain $\mathbb{E}_{\mathbb{W}_L}[\Phi_{L;u,v}(x,\cdot,\cdot)]$ by setting $f(t):=t^{2v},$ $\varphi(t):=\exp{\left(-e^{-2x}t\right)}.$ We first bound the innermost integral. Fix $D_2:=1/2.$ Then, using the bounds $|\sin(\cdot)|\leq 1, \sinh(p)\leq e^p/2, \cosh(p)\geq e^p/2,$ and  $-p^2 L^{-1}+p\leq L/4,$ we obtain
\begin{multline}\label{e:igb}
\left|\int_0^\infty \exp{\left(-\frac{p^2}{L}-yr\cosh(p)\right)}\sinh(p)\sin{\left(\frac{2\pi p}{L}\right)}dp\right|   
\leq\frac{1}{2}\int_0^\infty \exp{\left(-\frac{p^2}{L}+p-\frac{1}{2}yre^p\right)} dp
\\  \leq \frac{e^{L/4}}{2} \int_0^\infty \exp{\left(-D_2yre^p\right)} dp  = \frac{e^{L/4}}{2} \Gamma(0,D_2yr).
\end{multline}
Here, $\Gamma(0,z)$ denotes the upper incomplete gamma function defined for $z>0$ by 
\begin{align*}
\Gamma(0,z):=\int_z^\infty \frac{e^{-q}}{q}dq.
\end{align*}
Substituting this estimate into the preceding identity and absorbing the factor $e^{L/4}/2$ into $D_1$, we obtain
\begin{align}\label{e:YorGammaBound}
\mathbb{E}_{\mathbb{W}_L}\left[\Phi_{L;u,v}(x,\cdot,\cdot)\right] 
& \leq D_1 e^{-2(u+v)x} \int_0^\infty\int_0^\infty y^{2v} \exp{\left(-\frac{2e^{-2x}}{r}-\frac{r}{2}(1+y^2)\right)}\Gamma(0,D_2yr)drdy.
\end{align}

We first consider the case $v\geq0$. For $0<y\leq1$, we therefore have $y^{2v}\exp{\left(-\frac{r}{2}y^2\right)}\leq 1.$ Since, for every $z>0$, $\Gamma(0,z)\leq e^{-z}z^{-1},$ for $y\geq1$, we instead have $\Gamma(0,D_2ry)\leq e^{-D_2ry}(D_2ry)^{-1}.$ Splitting the $y$-integral, we obtain
\begin{multline}\label{e:b2}
\mathbb{E}_{\mathbb{W}_L}\left[\Phi_{L;u,v}(x,\cdot,\cdot)\right]  \leq D_1 e^{-2(u+v)x} \int_0^\infty \exp{\left( -\frac{2e^{-2x}}{r}-\frac{r}{2}\right)} \left( \int_0^1\Gamma(0,D_2ry) dy\right)dr
\\  +\frac{D_1}{D_2}e^{-2(u+v)x}\int_0^\infty \frac{1}{r} \exp{\left(-\frac{2e^{-2x}}{r}-\frac{r}{2}\right)}\left(\int_1^\infty y^{2v-1} \exp{\left(-\frac{r}{2}y^2-D_2ry\right)} dy\right) dr.
\end{multline}

We first bound the first term on the right-hand side of \eqref{e:b2}. Using the integral representation of the incomplete gamma function and Tonelli's theorem, we obtain, for every $r>0$,
\begin{align*}
\int_0^1\Gamma(0,D_2ry) dy & =\frac{1}{D_2r}\int_0^{D_2r}e^{-q} dq +\int_{D_2r}^\infty\frac{e^{-q}}{q} dq =\frac{1-e^{-D_2r}+D_2r\Gamma(0,D_2r)}{D_2r}.
\end{align*}
Since, for every $z>0$, $\Gamma(0,z)\leq e^{-z}z^{-1},$ we conclude that 
\begin{align*}
\int_0^1\Gamma(0,D_2ry) dy \leq \frac{1}{D_2r}.
\end{align*}
Consequently, the first term on the right-hand side of \eqref{e:b2} is bounded by
\begin{align*}
\frac{D_1}{D_2}e^{-2(u+v)x} \int_0^\infty \frac{1}{r} \exp{\left(-\frac{2e^{-2x}}{r}-\frac{r}{2}\right)}dr = \frac{2D_1}{D_2}e^{-2(u+v)x} K_0\left(2e^{-x}\right)
& \leq\frac{2D_1}{D_2}e^{-2(u+v)x}K_0\left(\sqrt2e^{-x}\right),
\end{align*}
where the inequality follows because $K_0$ is decreasing on $(0,\infty)$.

We next consider the second term on the right-hand side of \eqref{e:b2}. For $x\geq0$, it is simpler to bound the full expectation directly:
\begin{align*}
\mathbb{E}_{\mathbb{W}_L}\left[\Phi_{L;u,v}(x,\cdot,\cdot)\right]
& \leq e^{-2(u+v)x}\mathbb{E}_{\mathbb{W}_L}\left[e^{-2vg(L)}\right] =e^{Lv^2}e^{-2(u+v)x}.
\end{align*}
Since $K_0$ is positive and decreasing, for $x\geq0$, $K_0\left(\sqrt2e^{-x}\right) \geq K_0(\sqrt2)>0.$  Therefore, after possibly increasing the constant, the desired estimate follows for $x\geq0$.

It remains to consider $x<0$. For $v\geq0$ and $y\geq1$, $y^{2v-1}\leq y^{2v+1},$ and, therefore, for every $r>0$,
\begin{align*}
\int_1^\infty y^{2v-1}\exp{\left(-\frac{r}{2}y^2-D_2ry\right)}dy  & \leq \int_0^\infty y^{2v+1}\exp{\left(-\frac{r}{2}y^2\right)}dy = 2^v\Gamma(v+1)r^{-v-1}.
\end{align*}
Therefore, the second term on the right-hand side of \eqref{e:b2} is bounded by
\begin{align*}
\frac{2^vD_1\Gamma(v+1)}{D_2} e^{-2(u+v)x} \int_0^\infty r^{-v-2}\exp{\left(-\frac{2e^{-2x}}{r}-\frac{r}{2}\right)} dr.
\end{align*}

We set $A:=2e^{-2x}$ and then, since $x<0$, $A\geq2$. For $t>0$, there exists $C_v>0$ such that $t^{v+1}e^{-t}\leq C_ve^{-t/2}.$ It follows that $r^{-v-1}e^{-A/r} \leq C_ve^{-A/(2r)}$. Therefore, the second term on the right-hand side of \eqref{e:b2} is bounded by
\begin{align*}
C_{L,v}e^{-2(u+v)x} \int_0^\infty \frac{1}{r}\exp{\left(-\frac{e^{-2x}}{r}-\frac{r}{2}\right)}dr = 2C_{L,v}e^{-2(u+v)x}K_0\left(\sqrt{2}e^{-x}\right).
\end{align*}
Combining the preceding estimates, and possibly increasing the constant, we conclude that, when $v\geq0$,
\begin{align*}
\mathbb{E}_{\mathbb{W}_L}\left[\Phi_{L;u,v}(x,\cdot,\cdot)\right]\leq C_4e^{-2(u+v)x}K_0\left(\sqrt2e^{-x}\right).
\end{align*}

We now consider the case $v<0$. We split according to whether $e^{-g(L)}$ is smaller or larger than $1$:
\begin{align*}
\mathbb{E}_{\mathbb{W}_L}\left[\exp{\left(-2vg(L)-e^{-2x}\int_0^L e^{-2g(t)} dt\right)} \right]  & = \mathbb{E}_{\mathbb{W}_L}\left[\exp{\left(-2vg(L)-e^{-2x}\int_0^L e^{-2g(t)} dt\right)}\mathbf 1_{\{e^{-g(L)}<1\}}\right]
\\ & +\mathbb{E}_{\mathbb{W}_L}\left[\exp{\left(-2vg(L)-e^{-2x}\int_0^L e^{-2g(t)} dt\right)} \mathbf 1_{\{e^{-g(L)}\geq1\}} \right].
\end{align*}
We start with the second term. Applying the preceding Yor identity \eqref{e:YOR} with $f(y):=y^{2v}\mathbf 1_{\{y\geq1\}},$  $\varphi(t):=\exp{\left(-e^{-2x}t\right)},$ and using \eqref{e:igb}, we obtain
\begin{multline*}
\mathbb{E}_{\mathbb{W}_L}\left[\exp{\left(-2vg(L)-e^{-2x}\int_0^L e^{-2g(t)} dt\right)}\mathbf 1_{\{e^{-g(L)}\geq1\}}\right] 
\\ \leq \frac{D_1}{D_2}\int_0^\infty \frac{1}{r} \exp{\left(-\frac{2e^{-2x}}{r}-\frac{r}{2}\right)}\left(\int_1^\infty y^{2v-1} \exp{\left(-\frac{r}{2}y^2-D_2ry\right)}dy\right)dr.
\end{multline*}
Since $v<0$, $\int_1^\infty y^{2v-1}dy= (2|v|)^{-1}.$ Consequently,
\begin{align*}
\mathbb{E}_{\mathbb{W}_L}\left[\exp{\left(-2vg(L)-e^{-2x}\int_0^L e^{-2g(t)} dt\right)} \mathbf 1_{\{e^{-g(L)}\geq1\}}\right]
& \leq \frac{D_1}{2|v|D_2} \int_0^\infty\frac{1}{r}\exp{\left(-\frac{2e^{-2x}}r-\frac{r}{2}\right)}dr
\\ & =\frac{D_1}{|v|D_2}K_0\left(2e^{-x}\right) \leq \frac{D_1}{|v|D_2} K_0\left(\sqrt2e^{-x}\right).
\end{align*}

We next treat the first term. Since $e^{-g(L)}<1$ if and only if $g(L)>0,$ the Cameron--Martin formula gives
\begin{multline*}
\mathbb{E}_{\mathbb{W}_L}\left[\exp{\left(-2vg(L)-e^{-2x}\int_0^L e^{-2g(t)} dt\right)}\mathbf 1_{\{e^{-g(L)}<1\}}\right]
\\ =e^{Lv^2}\mathbb{E}_{\mathbb{W}_L}\left[\exp{\left(-e^{-2x}\int_0^L e^{-2(g(t)-vt)} dt\right)}\mathbf 1_{\{g(L)-vL>0\}}\right].
\end{multline*}
Because $v<0$, for every $t\in[0,L]$, $e^{2vt}\geq e^{2vL}.$ Therefore,
\begin{align*}
\int_0^L e^{-2(g(t)-vt)} dt =\int_0^L e^{-2g(t)}e^{2vt} dt \geq e^{2vL} \int_0^L e^{-2g(t)} dt.
\end{align*}
Dropping the indicator and using the monotonicity of the exponential, we obtain
\begin{align*}
\mathbb{E}_{\mathbb{W}_L}\left[\exp{\left(-2vg(L)-e^{-2x}\int_0^L e^{-2g(t)} dt\right)}\mathbf 1_{\{e^{-g(L)}<1\}}\right] 
& \leq e^{Lv^2}\mathbb{E}_{\mathbb{W}_L}\left[\exp{\left(-e^{-2(x-vL)}\int_0^L e^{-2g(t)} dt\right)}\right].
\end{align*}

The $v=0$ case proved above, applied with $u=0$ and with $x$ replaced by $x-vL$, implies that
\begin{align*}
\mathbb{E}_{\mathbb{W}_L}\left[\exp{\left(-2vg(L)-e^{-2x}\int_0^L e^{-2g(t)}dt\right)}\mathbf{1}_{\{e^{-g(L)}<1\}}\right] 
& \leq C_{L,v}K_0\left(\sqrt{2}e^{vL}e^{-x}\right).
\end{align*}

Combining the estimates on the two regions, we obtain
\begin{align}\label{e:vnegativeK0}
\mathbb{E}_{\mathbb{W}_L}\left[\exp{\left(-2vg(L)-e^{-2x}\int_0^L e^{-2g(t)} dt\right)}\right] 
& \leq C_{L,v}\left(K_0\left(\sqrt2e^{-x}\right)+K_0\left(\sqrt2e^{vL}e^{-x}\right)\right).
\end{align}

It remains to express the second term on the right-hand side in the form appearing in the statement of the lemma. We set $c:=\sqrt2e^{vL},$ and since $v<0$, $0<c<\sqrt2$. We claim that there exists $C_{L,v}>0$ such that, for every $z>0$,
\begin{align}\label{e:Besselcomparison}
K_0(cz)\leq C_{L,v}\left(K_0(\sqrt2z)+z^{|v|}K_{|v|}(cz)\right).
\end{align}
Indeed, as $z\to 0$, $K_0(cz)\sim-\log z,$ and $K_0(\sqrt2z)\sim-\log z,$ \cite[Entries 8.447.1, 8.447.3]{GR}
so the ratio of the first two quantities remains bounded near zero. On the other hand, as $z\to\infty$, $K_\nu(z)\sim\sqrt{\frac{\pi}{2z}}e^{-z}$ for every fixed $\nu$ \cite[ Entry 8.451.6]{GR}, and therefore $K_0(cz)(z^{|v|}K_{|v|}(cz))^{-1}\to 0.$ The inequality \eqref{e:Besselcomparison} then follows on the remaining compact range by continuity and positivity of the modified Bessel functions.

Applying \eqref{e:Besselcomparison} with $z=e^{-x}$, and using $K_{|v|}=K_v,$ and $(e^{-x})^{|v|}=e^{xv},$ we obtain
\begin{align*}
K_0\left(\sqrt{2}e^{vL}e^{-x}\right) \leq C_{L,v}\left(K_0\left(\sqrt2e^{-x}\right)+e^{xv}K_v\left(\sqrt2e^{vL}e^{-x}\right)\right).
\end{align*}
Combining this inequality with \eqref{e:vnegativeK0}, and multiplying by $e^{-2(u+v)x}$, we conclude that there exist constants $C_5,C_6>0$ such that
\begin{align*}
\mathbb{E}_{\mathbb{W}_L}\left[\Phi_{L;u,v}(x,\cdot,\cdot)\right] 
& \leq C_5 e^{-2(u+v)x}K_0\left(\sqrt2e^{-x}\right) +C_6 e^{-2(u+v)x+xv}K_v\left(C_7e^{-x}\right),
\end{align*}
where $C_7:=\sqrt{2}e^{vL}>0.$ This proves the claim when $v<0$ and completes the proof.
\end{proof}
As an immediate consequence of \cref{l:limitbounds}, we can show that $\mathcal{Z}_{L;u,v}$ is finite. In the case that $L=1$, this fact appears in~\cite{BKWW23}. This new proof is a novelty and new contribution of this paper. 
\begin{corollary}\label{c:part2}
Fix $L>0$ and $u,v\in\mathbb{R}$ such that $u+v>0$. Then $\mathcal Z_{L;u,v}<\infty.$
\end{corollary}
\begin{proof}
    When $v\geq 0$, 
    \begin{align*}
        \mathcal{Z}_{L;u,v}\leq C_{4}\int_{-\infty}^{\infty}e^{-2(u+v)x}K_{0}(\sqrt{2}e^{-x})dx.
    \end{align*} Identity $16$ in \cite{GR}, Section 6.561 states that for all $\mu,\nu,a$ such that $Re(\mu+1\pm\nu)>0$ and $Re(a)>0$,
\begin{align}\label{e:id3}
    \int_{0}^{\infty}r^{\mu}K_{\nu}(ar)dr = 2^{\mu-1}a^{-\mu-1}\Gamma\left(\frac{1+\mu+\nu}{2}\right)\Gamma\left(\frac{1+\mu-\nu}{2}\right).
\end{align} Making the change of coordinates $r=e^{-x}$, we can use this identity to evaluate the integral 
    \begin{align*}
       \mathcal{Z}_{L;u,v} \leq C_{4}\int_{0}^{\infty}r^{2(u+v)-1}K_{0}(\sqrt{2}r)dr =  C_{4}2^{(u+v)-2}\left(\Gamma\left(u+v\right)\right)^{2}.
    \end{align*}
    and when $v<0$, by \cref{l:limitbounds},
    \begin{align*}
        \mathcal{Z}_{L;u,v} & \leq C_{5}\int_{-\infty}^{\infty}e^{-2(u+v)x}K_{0}(\sqrt{2}e^{-x})dx  + C_{6}\int_{-\infty}^{\infty}e^{-2(u+v)x +xv}K_{v}\left(C_{7}e^{-x}\right)dx,
        \\ & \leq C_{5}2^{(u+v)-2}\left(\Gamma(u+v)\right)^{2}  + C_{6}2^{2(u+v)-v-2}C_{7}^{-(2u+v)}\Gamma\left(u+v\right)\Gamma\left(u\right).
    \end{align*} To evaluate the second integral we applied \eqref{e:id3}, the hypotheses of which are satisfied since $v<0$ and $u+v>0$ imply $u>0$. Thus both integrals are finite, and hence $\mathcal Z_{L;u,v}<\infty$. 
\end{proof} The preceding work also yields
\begin{align}\label{e:continuumlatticesum}
\sup_{N\geq1}N^{-1/2}\sum_{x\in\mathbb Z^{(N)}}\mathbb E_{\mathbb W_L}\left[\Phi_{L;u,v}(x,\cdot,\cdot)+\Phi_{L;u,v}(x+\log 2,\cdot,\cdot)\right]<\infty,
\end{align} which we record independently here for future use. 
\cref{c:part2} directly implies an additional corollary.
\begin{corollary}\label{c:finitemeas}
    Fix $L>0$ and $u,v\in\mathbb{R}$ such that $u+v>0$. Then the measure $\mathbb{Q}_{L;u,v}$ (\cref{d:kpzrn}) is a probability measure. 
\end{corollary} 
\section{Bounds on the Path Measure} \label{s:conv}

We want to show weak convergence of the sequence $\mathbb{Q}^{(N)}_{L;u,v}$. In this section, we prove bounds on the speed of convergence of the part of $\mathbb{Q}^{(N)}_{L;u,v}$ which depends on the path of the random walk. 
We consider an arbitrary bounded Lipschitz function $F(x,g,h)$, and set $C_{F}:=\sup_{x,g,h}|F(x,g,h)|+\mathrm{Lip}(F)=\|F\|_{\infty}+\mathrm{Lip}(F)$. This time, we use the standard KMT embedding, which helps us compare $\Phi^{(N)}_{L;u,v}(x,g,h)$ and $\Phi_{L;u,v}(x,g,h)$ on the same probability space. We note that 
\begin{multline}\label{e:bige} \big|\mathbb{E}_{\mathbb{W}^{(N)}_{L}}\left[F(x,\cdot,\cdot) \Phi^{(N)}_{L;u,v}(x,\cdot,\cdot)\right]-\mathbb{E}_{\mathbb{W}_{L}}\left[F(x,\cdot,\cdot) \Phi_{L;u,v}(x,\cdot,\cdot)\right] \big|
    \\  \leq C_{F}\mathbb{E}_{\mathbb{W}^{(N)}_{L}}\left[ \big|\Phi^{(N)}_{L;u,v}(x,\cdot,\cdot) -\Phi_{L;u,v}(x,\cdot,\cdot)\big|\right]  
     +  \Pi_{L;u,v}^{(N)}(x,F). 
\end{multline}
In the expression above, we define
\begin{align*}
    \Pi_{L;u,v}^{(N)}(x,F) & :=  \big| \mathbb{E}_{\mathbb{W}^{(N)}_{L}}\left[F(x,\cdot,\cdot)\Phi_{L;u,v}(x,\cdot,\cdot)\right] - \mathbb{E}_{\mathbb{W}_{L}}\left[F(x,\cdot,\cdot) \Phi_{L;u,v}(x,\cdot,\cdot)\right] \big|.
\end{align*} The subject of this section is to produce a bound on $\Pi_{L;u,v}^{(N)}(x,F)$. 
We use the KMT embedding theorem to show that this difference goes to $0$ very quickly as $N\to\infty$. In particular, we need to find a rate of convergence with a sufficiently nice dependence on $x$ so that we can use it to conclude  
\begin{align*}
    \big|\mathbb{E}_{\mathbb{P}^{(N)}_{L}}\left[F\cdot \Phi^{(N)}_{L;u,v}\right]-\mathbb{E}_{\mathbb{P}_{L}}\left[F\cdot \Phi_{L;u,v}\right]\big|\to 0.
\end{align*}
Now we state the adaptation of the KMT theorem which we will use.
\begin{theorem}[KMT coupling, adapted from \cite{KMT75}]
\label{t:KMT} Fix $L>0$. There exist constants $C_L,K_{L},\lambda_{L}>0$ such that, for every $N\in\mathbb{Z}_{>0}$ such that $\lfloor NL\rfloor\geq 1$, there is a probability space $(\Omega_{\textsc{KMT}}^{(N)},\mathcal F_{\textsc{KMT}}^{(N)},\mathbb{P}_{\textsc{KMT}}^{(N)})$ carrying $(C_0[0,L]\times C_0[0,L])$-valued random variables $\vec A^{(N)}=(A_1^{(N)},A_2^{(N)}),$ and $\vec B^{(N)}=(B_1^{(N)},B_2^{(N)})$ such that $\mathbb{P}_{\textsc{KMT}}^{(N)}\circ (\vec A^{(N)} )^{-1}=\mathbb{W}_L^{(N)}$ and $\mathbb{P}_{\textsc{KMT}}^{(N)}\circ(\vec B^{(N)})^{-1}=\mathbb{W}_L.$ Furthermore, for every $y\geq0$,
\begin{align*}
\mathbb{P}_{\textsc{KMT}}^{(N)}\left(\|\vec A^{(N)}-\vec B^{(N)}\|_{\infty,2}\geq C_L N^{-1/2}(\log N+y)\right)\leq K_{L}e^{-\lambda_{L}y}.
\end{align*}
\end{theorem}
With this theorem in hand, we can prove the following bound:
\begin{proposition}\label{p:pin} Fix $L>0$ and $u,v\in\mathbb{R}$ such that $u+v>0$. 
     There exists $N_{0}:=N_{0}(L,u,v)\in\mathbb{Z}_{>0}$ such that for all bounded Lipschitz functions $F(x,g,h)$ there exists $C:=C(L,u,v,C_{F})>0$ for which, for every $N\geq N_{0}$  and $x\in\mathbb{R}$, 
    \begin{align*}
        \Pi^{(N)}_{L;u,v}(x,F) \leq  C\log{(N)}N^{-1/2}\left(\mathbb{E}_{\mathbb{W}_{L}}\left[\Phi_{L;u,v}(x,\cdot,\cdot)+\Phi_{L;u,v}(x+\log{(2)},\cdot,\cdot)\right]\right) + Ce^{-2(u+v)x}N^{-(u+v+1)}. 
    \end{align*}
\end{proposition}
\begin{proof} 
For all fixed $N\in\mathbb{Z}_{>0}$ such that $\lfloor LN\rfloor\geq 1$, let $(\Omega_{\textsc{KMT}}^{(N)},
\mathcal F_{\textsc{KMT}}^{(N)},
\mathbb{P}_{\textsc{KMT}}^{(N)})$, $\vec A^{(N)}=(A_1^{(N)},A_2^{(N)}),$ $\vec B^{(N)}=(B_1^{(N)},B_2^{(N)})$, and $C_{L},K_{L},\lambda_{L}>0$ be as in \cref{t:KMT}. We will use the notation $\mathbb{E}_{\textsc{KMT}}^{(N)}$ to denote expectation with respect to $\mathbb{P}_{\textsc{KMT}}^{(N)}$. 

For fixed $x\in\mathbb{R}$ and for $y\geq 0$, we define the event $$K_{N,y}:=\left\{\|\vec A^{(N)}-\vec B^{(N)}\|_{\infty,2}\geq C_{L}N^{-\frac{1}{2}}(\log{(N)}+y)\right\}.$$
Then we can bound $
\Pi_{L;u,v}^{(N)}(x,F)$ in terms of $K_{N,y}$ and its complement. We find that $
\Pi_{L;u,v}^{(N)}(x,F)$ is bounded above by
\begin{multline}\label{e:kmt} 
\leq  \mathbb{E}_{\textsc{KMT}}^{(N)}\left[\big|F(x,A_{1}^{(N)},A_{2}^{(N)})\Phi_{L;u,v}(x,A_{1}^{(N)},A_{2}^{(N)})-F(x,B_{1}^{(N)},B_{2}^{(N)})\Phi_{L;u,v}(x,B_{1}^{(N)},B_{2}^{(N)})\big|\mathbf{1}_{K_{N,y}^{c} }\right] \\  + \mathbb{E}_{\textsc{KMT}}^{(N)}\left[\big|F(x,A_{1}^{(N)},A_{2}^{(N)})\Phi_{L;u,v}(x,A_{1}^{(N)},A_{2}^{(N)})-F(x,B_{1}^{(N)},B_{2}^{(N)})\Phi_{L;u,v}(x,B_{1}^{(N)},B_{2}^{(N)})\big|\mathbf{1}_{K_{N,y}}\right]
    \end{multline}
    We first bound the term associated with the event $K_{N,y}$. Since  
    \begin{align*}
        |F(x,g,h)\Phi_{L;u,v}(x,g,h)|\leq C_{F}e^{-2(u+v)x-2vg(L)},
    \end{align*} then for all $r>0$, we can write 
    \begin{multline*}
        \mathbb{E}_{\textsc{KMT}}^{(N)}\left[\big|F(x,A_{1}^{(N)},A_{2}^{(N)})\Phi_{L;u,v}(x,A_{1}^{(N)},A_{2}^{(N)})-F(x,B_{1}^{(N)},B_{2}^{(N)})\Phi_{L;u,v}(x,B_{1}^{(N)},B_{2}^{(N)})\big|\mathbf{1}_{K_{N,y}}\right]
        \\ \leq 2C_{F}e^{-2(u+v)x+2|v|r}\mathbb{P}_{\textsc{KMT}}^{(N)}(K_{N,y}) + C_{F}e^{-2(u+v)x}\mathbb{E}_{\textsc{KMT}}^{(N)}\left[e^{-2vA_{1}^{(N)}(L)}\mathbf{1}_{|A_{1}^{(N)}(L)|>r}\right]
        \\ + C_{F}e^{-2(u+v)x}\mathbb{E}_{\textsc{KMT}}^{(N)}\left[e^{-2vB_{1}^{(N)}(L)}\mathbf{1}_{|B_{1}^{(N)}(L)|>r}\right].
    \end{multline*} Furthermore, for any $v\in\mathbb{R}$, 
    \begin{align*}
        \mathbb{E}_{\textsc{KMT}}^{(N)}\left[e^{-2vA_{1}^{(N)}(L)}\right] = \left(\frac{1+\cosh{(-2vN^{-1/2})}}{2}\right)^{\lfloor LN\rfloor}\leq e^{Lv^{2}}, & & \mathbb{E}_{\textsc{KMT}}^{(N)}\left[e^{-2vB_{1}^{(N)}(L)}\right] = e^{Lv^{2}}.
    \end{align*} Therefore, we can apply a Chernoff bound to obtain for all $r\geq 4L|v|,$
    \begin{align*}
        \mathbb{E}_{\textsc{KMT}}^{(N)}\left[e^{-2vA_{1}^{(N)}(L)}\mathbf{1}_{|A_{1}^{(N)}(L)|>r}\right], \mathbb{E}_{\textsc{KMT}}^{(N)}\left[e^{-2vB_{1}^{(N)}(L)}\mathbf{1}_{|B_{1}^{(N)}(L)|>r}\right] \leq 2e^{-r^{2}/2L}.
    \end{align*} This implies that 
    \begin{multline*}
        \mathbb{E}_{\textsc{KMT}}^{(N)}\left[\big|F(x,A_{1}^{(N)},A_{2}^{(N)})\Phi_{L;u,v}(x,A_{1}^{(N)},A_{2}^{(N)})-F(x,B_{1}^{(N)},B_{2}^{(N)})\Phi_{L;u,v}(x,B_{1}^{(N)},B_{2}^{(N)})\big|\mathbf{1}_{K_{N,y}}\right]
        \\ \leq 2C_{F}K_{L}e^{-2(u+v)x+2|v|r-\lambda_{L}y} + 4C_{F}e^{-2(u+v)x-r^{2}/2L}.
    \end{multline*} Then, setting $r=(\lambda_{L}Ly)^{1/2}$, there exists a constant $C_{F,L,v}>0$ such that both terms on the right hand side are bounded by $C_{F,L,v}e^{-2(u+v)x-\lambda_{L}y/2}.$ We then choose $y:=y(N)=2(u+v+1)\lambda_{L}^{-1}\log{(N)}$, which results in a bound of 
    \begin{multline}\label{e:secondterm}
        \mathbb{E}_{\textsc{KMT}}^{(N)}\left[\big|F(x,A_{1}^{(N)},A_{2}^{(N)})\Phi_{L;u,v}(x,A_{1}^{(N)},A_{2}^{(N)})-F(x,B_{1}^{(N)},B_{2}^{(N)})\Phi_{L;u,v}(x,B_{1}^{(N)},B_{2}^{(N)})\big|\mathbf{1}_{K_{N,y}}\right]
        \\ \leq C_{F,L,v}e^{-2(u+v)x}N^{-(u+v+1)}. 
    \end{multline} 
Next, we bound the term in \eqref{e:kmt} associated with $K_{N,y}^{c}$. For some $\delta>0$, let $g,h,\widetilde g,\widetilde h\in C_0[0,L]$ satisfy $\|(g,h)-(\widetilde g,\widetilde h)\|_{\infty,2}<\delta.$
Since $F$ is Lipschitz, this implies that $|F(x,g,h)-F(x,\widetilde g,\widetilde h)|
\leq\mathrm{Lip}(F)\delta.$
Consequently,
\begin{multline*}
| F(x,g,h)\Phi_{L;u,v}(x,g,h) - F(x,\widetilde{g},\widetilde{h}) \Phi_{L;u,v}(x,\widetilde{g},\widetilde{h}) | 
\\ \leq C_{F}\delta \Phi_{L;u,v}(x,\widetilde g,\widetilde h) + C_{F}| \Phi_{L;u,v}(x,g,h) - \Phi_{L;u,v}(x,\widetilde g,\widetilde h) |.
\end{multline*}
In particular, since $\|g-\widetilde g\|_\infty<\delta$ and because $\Phi_{L;u,v}$ does not depend on its third argument, we conclude that 
\begin{multline}\label{e:ftotbounds}
|\Phi_{L;u,v}(x,g,h)-\Phi_{L;u,v}(x,\widetilde g,\widetilde h)|
\\ = \Phi_{L;u,v}(x,\widetilde g,\widetilde h)\left|\exp{\left(-2v(g(L)-\widetilde g(L))-e^{-2x}\int_0^Le^{-2\widetilde g(t)}\left(e^{-2(g(t)-\widetilde g(t))}-1\right)dt\right)}-1\right|.
\end{multline}
We define $m:=e^{-2x}\int_0^L e^{-2\widetilde g(t)}dt.$ Since $|g(L)-\widetilde g(L)|<\delta$ and $e^{-2\delta}-1\leq e^{-2(g(t)-\widetilde g(t))}-1\leq e^{2\delta}-1,$ the final absolute value in \eqref{e:ftotbounds} is bounded by
\begin{align*}G_{\delta}(m):=\max\left\{\exp{\left(2|v|\delta+\left(1-e^{-2\delta}\right)m\right)}-1,1-\exp{\left(-2|v|\delta-\left(e^{2\delta}-1\right)m\right)}\right\}.
\end{align*}
Therefore, we conclude that 
\begin{align*}
|\Phi_{L;u,v}(x,g,h)-\Phi_{L;u,v}(x,\widetilde g,\widetilde h) |\leq \Phi_{L;u,v}(x,\widetilde g,\widetilde h) G_{\delta}\left(e^{-2x}\int_0^L e^{-2\widetilde g(t)}dt\right).
\end{align*}
We note that there exist $\delta_0,D_v>0$ such that for every $0<\delta\leq\delta_0$ and $m\geq0$,
\begin{align}\label{e:Gdelta}
e^{-m}G_\delta(m) \leq D_v\delta e^{-m/4}.
\end{align}
We justify this assertion by noting that for all sufficiently small $\delta$, each of the following inequalities must hold
\begin{align*}
e^{-2\delta}\geq 3/4, & &  1-e^{-2\delta}\leq 2\delta,  & &  e^{2\delta}-1\leq 4\delta.
\end{align*}
Using the inequalities $e^z-1\leq ze^z$ and $1-e^{-z}\leq z$ for $z\geq 0$, we conclude that the two terms inside the maximum which defines $G_\delta(m)$, after multiplication by $e^{-m}$, are bounded respectively by $C_v\delta(1+m)e^{-3m/4}$ and $C_v\delta(1+m)e^{-m},$ and both expressions are bounded by $D_v\delta e^{-m/4}$.

Therefore, since
\begin{align*}
    \Phi_{L;u,v}(x,\widetilde g,\widetilde h)=e^{-2(u+v)x-2v\widetilde g(L)}e^{-m},
\end{align*}
we can apply \eqref{e:Gdelta} to see 
\begin{align*}
|\Phi_{L;u,v}(x,g,h)- \Phi_{L;u,v}(x,\widetilde g,\widetilde h)|\leq D_v\delta e^{-2(u+v)x-2v\widetilde g(L)-m/4}.
\end{align*}
Furthermore, since $\Phi_{L;u,v}\left(x+\log{(2)},\widetilde{g},\widetilde{h}\right)= 2^{-2(u+v)}e^{-2(u+v)x-2v\widetilde g(L)-m/4}$, there exists $D_{u,v}>0$ such that, for every $0<\delta\leq\delta_0$,
\begin{align}\label{e:phipathcomparison}
\left|\Phi_{L;u,v}(x,g,h)-\Phi_{L;u,v}(x,\widetilde g,\widetilde h)\right|  & \leq D_{u,v}\delta\Phi_{L;u,v}\left(x+\log{(2)},\widetilde{g},\widetilde{h}\right).
\end{align}

We now apply the preceding deterministic estimate with $(g,h)= (A_1^{(N)},A_2^{(N)})$, $(\widetilde g,\widetilde h)=(B_1^{(N)},B_2^{(N)})$, $\delta:=\delta(N)= C_LN^{-1/2}(\log{(N)}+y)$, and, as before, $y:=y(N)=2(u+v+1)\lambda_L^{-1}\log{(N)},$ we have $\delta\leq C_{L,u,v}\log{(N)}N^{-1/2},$ and for all $N$ sufficiently large, $\delta\leq\delta_0$. On $K_{N,y }^{c}$, the hypotheses of
\eqref{e:phipathcomparison} are therefore satisfied. Therefore, on $K_{N,y}^{c}$,
\begin{multline*}
 |F(x,A_1^{(N)},A_2^{(N)})\Phi_{L;u,v}(x,A_1^{(N)},A_2^{(N)}) -F(x,B_1^{(N)},B_2^{(N)})\Phi_{L;u,v}(x,B_1^{(N)},B_2^{(N)})|
\\ \leq C_F\delta\Phi_{L;u,v}\left(x,B_1^{(N)},B_2^{(N)}\right) + C_F D_{u,v}\delta \Phi_{L;u,v}\left(x+\log(2),B_1^{(N)},B_2^{(N)}\right).
\end{multline*}
Taking expectations, we obtain
\begin{multline*}
\mathbb{E}_{\textsc{KMT}}^{(N)}\left[|F(x,A_1^{(N)},A_2^{(N)})\Phi_{L;u,v}(x,A_1^{(N)},A_2^{(N)}) -F(x,B_1^{(N)},B_2^{(N)})\Phi_{L;u,v}(x,B_1^{(N)},B_2^{(N)})|\mathbf{1}_{K_{N,y}^{c}}\right]
\\ \leq C_{L,u,v,C_F}\delta \mathbb{E}_{\mathbb{W}_L}\left[\Phi_{L;u,v}(x,\cdot,\cdot)+\Phi_{L;u,v}(x+\log(2),\cdot,\cdot)\right].
\end{multline*}
Then, using the bound $\delta \leq C_{L,u,v}\log{ (N)}N^{-1/2},$ we conclude that
\begin{multline}\label{e:firstterm}
 \mathbb{E}_{\textsc{KMT}}^{(N)}\left[|F(x,A_1^{(N)},A_2^{(N)})\Phi_{L;u,v}(x,A_1^{(N)},A_2^{(N)}) -F(x,B_1^{(N)},B_2^{(N)})\Phi_{L;u,v}(x,B_1^{(N)},B_2^{(N)})|\mathbf{1}_{K_{N,y}^{c}}\right]
\\ \leq C_{L,u,v,C_F}\log{(N)}N^{-1/2}\mathbb{E}_{\mathbb{W}_L}\left[\Phi_{L;u,v}(x,\cdot,\cdot)+\Phi_{L;u,v}(x+\log(2),\cdot,\cdot)\right].
\end{multline}
We set $C:=\max\left\{C_{L,v,C_F},C_{L,u,v,C_F}\right\}$ and choose $N_0=N_0(L,u,v)$ sufficiently large that all of the preceding estimates hold. Then \eqref{e:kmt}, \eqref{e:firstterm}, and \eqref{e:secondterm} imply
\begin{align*}
\Pi_{L;u,v}^{(N)}(x,F) & \leq  C \log {(N)}N^{-1/2} \mathbb{E}_{\mathbb{W}_L}\left[\Phi_{L;u,v}(x,\cdot,\cdot) +\Phi_{L;u,v}(x+\log 2,\cdot,\cdot) \right] + Ce^{-2(u+v)x}N^{-(u+v+1)}.
\end{align*}
This proves \cref{p:pin}.
\end{proof}
\section{Weak Convergence}\label{s:weakconvergencefinal}
Throughout this section, fix $L>0$ and $u,v\in\mathbb{R}$ such that
$u+v>0$. We prove \cref{t:main}, beginning with the following
intermediate result.

\begin{proposition}\label{p:almost} For all bounded Lipschitz functions $F:\mathbb{R}\times C_{0}[0,L]\times C_{0}[0,L]\to\mathbb{R},$
\begin{align*}
    \lim_{N\to\infty}\mathbb{E}_{\mathbb{P}^{(N)}_{L}}[F\cdot \Phi^{(N)}_{L;u,v}] = \mathbb{E}_{\mathbb{P}_{L}}[F\cdot \Phi_{L;u,v}].
\end{align*}
\end{proposition}
Before proving \cref{p:almost}, we give the proof of \cref{t:main} using this result. We state the following corollary of \cref{p:almost}, which we obtain by setting $F\equiv 1$. 
\begin{corollary}\label{c:partition}
    The partition functions converge, $\lim_{N\to\infty}\mathcal{Z}^{(N)}_{L;u,v}=\mathcal{Z}_{L;u,v}$.
\end{corollary} With these propositions, we can prove \cref{t:main}.
\begin{proof}[Proof of \cref{t:main}]
We combine the results of \cref{c:partition} and \cref{p:almost}. Since $\mathcal{Z}_{L;u,v}>0$, this proves that for all bounded Lipschitz functions $F:\mathbb{R}\times C_{0}[0,L]\times C_{0}[0,L]\to\mathbb{R},$ $$\lim_{N\to\infty}(\mathcal{Z}^{(N)}_{L;u,v})^{-1}\mathbb{E}_{\mathbb{P}^{(N)}_{L}}\left[F\cdot \Phi^{(N)}_{L;u,v}\right]=(\mathcal{Z}_{L;u,v})^{-1}\mathbb{E}_{\mathbb{P}_{L}}\left[F\cdot \Phi_{L;u,v}\right].$$
Since bounded Lipschitz functions determine weak convergence on
$\mathbb{R}\times C_{0}[0,L]\times C_{0}[0,L]$, this proves
\cref{t:main}.
\end{proof}
Finally, we finish the proof of \cref{p:almost}.
\begin{proof}[Proof of \cref{p:almost}] 
We use $x$ to denote lattice points in $\mathbb Z^{(N)}$ and $r$ to denote elements of $\mathbb R$. We denote the minimal element of $\mathbb Z^{(N)}$ by $\widetilde{x}_N:=N^{-1/2}-\log(\sqrt N).$ 
Expanding the expression for the difference and then applying the triangle inequality, we obtain
\begin{multline}\label{e:firstbound}
\left|\mathbb E_{\mathbb P_L^{(N)}}\left[F\cdot\Phi_{L;u,v}^{(N)}\right]-\mathbb E_{\mathbb P_L}\left[F\cdot\Phi_{L;u,v}\right]\right|
\\ \leq N^{-1/2}\sum_{x\in\mathbb Z^{(N)}}\bigg|\mathbb E_{\mathbb W_L^{(N)}}\left[F(x,\cdot,\cdot)\Phi_{L;u,v}^{(N)}(x,\cdot,\cdot)\right]-\mathbb E_{\mathbb W_L}\left[F(x,\cdot,\cdot)\Phi_{L;u,v}(x,\cdot,\cdot)\right]\bigg| 
\\ +\sum_{x\in\mathbb Z^{(N)}}\int_x^{x+N^{-1/2}}\bigg|\mathbb E_{\mathbb W_L}\left[F(x,\cdot,\cdot)\Phi_{L;u,v}(x,\cdot,\cdot)\right]-\mathbb E_{\mathbb W_L}\left[F(r,\cdot,\cdot)\Phi_{L;u,v}(r,\cdot,\cdot)\right]\bigg|dr
\\ +\int_{-\infty}^{\widetilde{x}_N}\bigg|\mathbb E_{\mathbb W_L}\left[F(r,\cdot,\cdot)\Phi_{L;u,v}(r,\cdot,\cdot)\right]\bigg|dr.
\end{multline}
Since $\Phi_{L;u,v}$ is nonnegative, the third term on the right-hand side of \eqref{e:firstbound} is bounded as
\begin{align*}
\int_{-\infty}^{\widetilde{x}_N}\bigg|\mathbb E_{\mathbb W_L}\left[F(r,\cdot,\cdot)\Phi_{L;u,v}(r,\cdot,\cdot)\right]\bigg|dr
\leq \|F\|_{\infty}\int_{-\infty}^{\widetilde{x}_N}\mathbb E_{\mathbb W_L}\left[\Phi_{L;u,v}(r,\cdot,\cdot)\right]dr.
\end{align*}
By \cref{c:part2}, $\mathcal Z_{L;u,v}<\infty.$ Since $\widetilde{x}_N\to-\infty$ as $N\to\infty$, it therefore follows that
\begin{align}\label{e:lowertail}
\lim_{N\to\infty}\int_{-\infty}^{\widetilde{x}_N}\bigg|\mathbb E_{\mathbb W_L}\left[F(r,\cdot,\cdot)\Phi_{L;u,v}(r,\cdot,\cdot)\right]\bigg|dr=0.
\end{align}
By \eqref{e:bige}, the first term on the right-hand side of \eqref{e:firstbound} is bounded by
\begin{align}\label{e:pathmeasureterms}
\|F\|_{\infty}N^{-1/2}\sum_{x\in\mathbb Z^{(N)}}\mathbb E_{\mathbb W_L^{(N)}}\left[\left|\Phi_{L;u,v}^{(N)}(x,\cdot,\cdot)-\Phi_{L;u,v}(x,\cdot,\cdot)\right|\right]  +N^{-1/2}\sum_{x\in\mathbb Z^{(N)}}\Pi_{L;u,v}^{(N)}(x,F).
\end{align}
Then, by \cref{p:pin}, the second term of \eqref{e:pathmeasureterms} is bounded by
\begin{align*}
C\log(N)N^{-1}\sum_{x\in\mathbb Z^{(N)}}\mathbb E_{\mathbb W_L}\left[\Phi_{L;u,v}(x,\cdot,\cdot)+\Phi_{L;u,v}(x+\log 2,\cdot,\cdot)\right]+CN^{-1/2-(u+v+1)}\sum_{x\in\mathbb Z^{(N)}}e^{-2(u+v)x}.
\end{align*}
By \eqref{e:continuumlatticesum}, the first term in the preceding
bound is $O(\log(N)N^{-1/2})$ and therefore converges to zero.

Furthermore,
\begin{align*}
N^{-(u+v+1)}N^{-1/2}\sum_{x\in\mathbb Z^{(N)}}e^{-2(u+v)x} =N^{-3/2}\sum_{j=1}^{\infty}e^{-2(u+v)j/\sqrt N} =O(N^{-1}),
\end{align*}
and, therefore, the second term in \eqref{e:pathmeasureterms} goes to $0$ as $N\to\infty$.

It remains to show that the first term of \eqref{e:pathmeasureterms} also goes to $0$ as $N\to\infty$. We fix $R>0$ and claim that
\begin{align}\label{e:weightdifferencecompact}
\lim_{N\to\infty} \sup_{x\in\mathbb Z^{(N)}\cap[-R,R]}\mathbb E_{\mathbb W_L^{(N)}}\left[\left|\Phi_{L;u,v}^{(N)}(x,\cdot,\cdot)-\Phi_{L;u,v}(x,\cdot,\cdot)\right|\right] =0.
\end{align}
By way of contradiction, suppose that \eqref{e:weightdifferencecompact} does not hold. Then there exists $\varepsilon>0$, as well as sequences $N_k\to\infty$ and $x_k\in\mathbb Z^{(N_k)}\cap[-R,R]$ such that for all $k\geq 1$,
\begin{align}\label{e:compactlowerbound}
\mathbb E_{\mathbb W_L^{(N_k)}}\left[\left|\Phi_{L;u,v}^{(N_k)}(x_k,\cdot,\cdot)-\Phi_{L;u,v}(x_k,\cdot,\cdot)\right|\right]\geq\varepsilon.
\end{align}
Possibly passing to a subsequence, we assume that $x_k$ converges to $x\in[-R,R]$.
Since $\mathbb W_L^{(N_k)}$ converges weakly to $\mathbb W_L$, we apply the Skorokhod representation theorem \cite[Theorem 6.7]{Bill99} to obtain a joint probability space $(\Omega_{\mathrm{Sk}},\mathcal F_{\mathrm{Sk}},
\mathbb P_{\mathrm{Sk}})$ on which  $(g_k,h_k)$ have law $\mathbb W_L^{(N_k)}$ and $(g,h)$ have law $\mathbb W_L$ and such that $(g_k,h_k)$ converges uniformly to $(g,h)$ on $[0,L]$ almost surely in $\mathbb P_{\mathrm{Sk}}$. On this probability space,
\begin{align*}
\lim_{k\to\infty} \max_{0\leq i\leq\lfloor LN_k\rfloor}\frac{e^{-2(g_k(t_i)+x_k)}}{N_k}= 0, & &  \lim_{k\to\infty} \frac{1}{N_k} \sum_{i=0}^{\lfloor LN_k\rfloor} e^{-2(g_k(t_i)+x_k)} = e^{-2x}\int_0^L e^{-2g(t)}dt,
\end{align*}
\begin{align*}
     \lim_{k\to\infty} \frac{1}{N_k^2} \sum_{i=0}^{\lfloor LN_k\rfloor} e^{-4(g_k(t_i)+x_k)}= 0,
\end{align*}
almost surely in $\mathbb P_{\mathrm{Sk}}$. Therefore, using $\log(1-z)=-z+O(z^2)$, we obtain
\begin{align*}
\lim_{k\to\infty} \prod_{i=0}^{\lfloor LN_k\rfloor}\left(1-\frac{e^{-2(g_k(t_i)+x_k)}}{N_k}\right) = \exp\left(-e^{-2x}\int_0^L e^{-2g(t)}dt\right)
\end{align*}
almost surely in $\mathbb P_{\mathrm{Sk}}$. Finally, this implies that almost surely in $\mathbb P_{\mathrm{Sk}}$,
\begin{align*}
\lim_{k\to\infty} \left|\Phi_{L;u,v}^{(N_k)}(x_k,g_k,h_k)-\Phi_{L;u,v}(x_k,g_k,h_k)\right| = 0.
\end{align*}
For all $x_k\in[-R,R]$,
\begin{align*}
\Phi_{L;u,v}^{(N_k)}(x_k,g_k,h_k), \Phi_{L;u,v}(x_k,g_k,h_k)\leq \exp\left(2(u+v)R+2|v||g_k(L)|\right).
\end{align*}
The endpoint variables $g_k(L)$ have exponential moments of every fixed order that are bounded uniformly in $k$. Hence the preceding sequence is uniformly integrable. This justifies the convergence in expectation,
\begin{align*}
    \lim_{k\to\infty} \mathbb E_{\mathbb W_L^{(N_k)}}\left[\left|\Phi_{L;u,v}^{(N_k)}(x_k,\cdot,\cdot)-\Phi_{L;u,v}(x_k,\cdot,\cdot)\right|\right] = 0.
\end{align*}
This is in contradiction with \eqref{e:compactlowerbound} and therefore proves
\eqref{e:weightdifferencecompact}.
Since $|\mathbb{Z}^{(N)}\cap[-R,R]|\leq 2R\sqrt{N}+1$,  \eqref{e:weightdifferencecompact} implies
\begin{align}\label{e:weightdifferencecompactsum}
\lim_{N\to\infty} N^{-1/2} \sum_{x\in\mathbb Z^{(N)}\cap[-R,R]} \mathbb E_{\mathbb W_L^{(N)}}\left[\left|\Phi_{L;u,v}^{(N)}(x,\cdot,\cdot)-\Phi_{L;u,v}(x,\cdot,\cdot)\right|\right] = 0.
\end{align}

What remains is to apply \eqref{e:weightdifferencecompactsum} to prove  
\begin{align}\label{e:weightdifference}
\lim_{N\to\infty} N^{-1/2}\sum_{x\in\mathbb Z^{(N)}}\mathbb E_{\mathbb W_L^{(N)}}\left[\left|\Phi_{L;u,v}^{(N)}(x,\cdot,\cdot)-\Phi_{L;u,v}(x,\cdot,\cdot)\right|\right] =0.
\end{align}
We fix $R>0$, and note that
\begin{multline*}
N^{-1/2}\sum_{x\in\mathbb Z^{(N)}}\mathbb E_{\mathbb W_L^{(N)}}\left[\left|\Phi_{L;u,v}^{(N)}(x,\cdot,\cdot)-\Phi_{L;u,v}(x,\cdot,\cdot)\right|\right]\leq N^{-1/2}\sum_{\substack{x\in\mathbb Z^{(N)}\\ |x|\leq R}}\mathbb E_{\mathbb W_L^{(N)}}\left[\left|\Phi_{L;u,v}^{(N)}(x,\cdot,\cdot)-\Phi_{L;u,v}(x,\cdot,\cdot)\right|\right] 
\\ + N^{-1/2} \sum_{\substack{x\in\mathbb Z^{(N)}\\ |x|>R}}\mathbb E_{\mathbb W_L^{(N)}}\left[\Phi_{L;u,v}^{(N)}(x,\cdot,\cdot)\right] + N^{-1/2}\sum_{\substack{x\in\mathbb Z^{(N)}\\ |x|>R}}\mathbb E_{\mathbb W_L^{(N)}}\left[\Phi_{L;u,v}(x,\cdot,\cdot)\right].
\end{multline*}
For each fixed $R$, the first term on the right-hand side converges to zero by \eqref{e:weightdifferencecompact} and the bound $N^{-1/2}|\mathbb Z^{(N)}\cap[-R,R]|\leq 2R+N^{-1/2}$. By \cref{p:pin}, applied with $F\equiv1$,
\begin{align*}
N^{-1/2}\sum_{\substack{x\in\mathbb Z^{(N)}\\ |x|>R}}\mathbb E_{\mathbb W_L^{(N)}}\left[\Phi_{L;u,v}(x,\cdot,\cdot)\right] 
& \leq N^{-1/2}\sum_{\substack{x\in\mathbb Z^{(N)}\\ |x|>R}}\mathbb E_{\mathbb W_L}\left[\Phi_{L;u,v}(x,\cdot,\cdot)\right]+N^{-1/2}\sum_{\substack{x\in\mathbb Z^{(N)}\\ |x|>R}}\Pi_{L;u,v}^{(N)}(x,1).
\end{align*}
The final term is bounded above by $N^{-1/2}\sum_{x\in\mathbb Z^{(N)}}\Pi_{L;u,v}^{(N)}(x,1).$ By \cref{p:pin}, \eqref{e:continuumlatticesum}, and the preceding calculation, this expression is $O(\log(N)N^{-1/2}+N^{-1})$, and hence converges to zero.
Consequently,
\begin{multline*}
\limsup_{N\to\infty}N^{-1/2}\sum_{x\in\mathbb Z^{(N)}}\mathbb E_{\mathbb W_L^{(N)}}\left[\left|\Phi_{L;u,v}^{(N)}(x,\cdot,\cdot)-\Phi_{L;u,v}(x,\cdot,\cdot)\right|\right]
\\ \leq\limsup_{N\to\infty}N^{-1/2}\sum_{\substack{x\in\mathbb Z^{(N)}\\ |x|>R}}\mathbb E_{\mathbb W_L^{(N)}}\left[\Phi_{L;u,v}^{(N)}(x,\cdot,\cdot)\right] +\limsup_{N\to\infty}N^{-1/2}\sum_{\substack{x\in\mathbb Z^{(N)}\\ |x|>R}}\mathbb E_{\mathbb W_L}\left[\Phi_{L;u,v}(x,\cdot,\cdot)\right].
\end{multline*} 
The estimates in \cref{l:finitebounds,l:limitbounds} bound the summands by integrable functions of $x$ that are independent of $N$, apart from the $N$-dependent remainder in \cref{l:finitebounds}, whose entire normalized lattice sum is $O(N^{-1})$ by the preceding calculation. Applying \eqref{e:continuumlatticesum} to the sums restricted to
$|x|>R$ shows that both terms on the right-hand side converge to zero as $R\to\infty$. Thus, letting $R\to\infty$ proves \eqref{e:weightdifference}.
Together with the preceding estimate for the second term in \eqref{e:pathmeasureterms}, this shows that the first term on the right-hand side of \eqref{e:firstbound} converges to zero.

We finally consider the second term on the right-hand side of \eqref{e:firstbound}. For $x\in\mathbb Z^{(N)}$ and $r\in[x,x+N^{-1/2}]$, we have
\begin{align*}
\bigg|\mathbb E_{\mathbb W_L}\left[F(x,\cdot,\cdot)\Phi_{L;u,v}(x,\cdot,\cdot)\right]-\mathbb E_{\mathbb W_L}\left[F(r,\cdot,\cdot)\Phi_{L;u,v}(r,\cdot,\cdot)\right]\bigg| 
&\leq\mathrm{Lip}(F)(r-x) \mathbb E_{\mathbb W_L}\left[\Phi_{L;u,v}(x,\cdot,\cdot)\right] 
\\ & + \|F\|_{\infty}\mathbb E_{\mathbb W_L}\left[\left|\Phi_{L;u,v}(x,\cdot,\cdot)-\Phi_{L;u,v}(r,\cdot,\cdot)\right|\right].
\end{align*}
For $r\in[x,x+N^{-1/2}]$, we instead compare the two weights directly. There exists $C>0$, independent of $N,x,r,g,$ and $h$, such that
\begin{multline*}
\left|\Phi_{L;u,v}(x,g,h)-\Phi_{L;u,v}(r,g,h)\right| \leq \left|e^{-2(u+v)x}-e^{-2(u+v)r}\right|e^{-2vg(L)}\exp\left(-e^{-2x}\int_0^L e^{-2g(t)}dt\right) \\ +e^{-2(u+v)r-2vg(L)}\left|\exp\left(-e^{-2x}\int_0^L e^{-2g(t)}dt\right)-\exp\left(-e^{-2r}\int_0^L e^{-2g(t)}dt\right)\right|.
\end{multline*}
Since $e^{2(u+v)(r-x)}-1\leq C(r-x)$ and 
\begin{align*}\exp\left(-e^{-2x}\int_0^L e^{-2g(t)}dt\right)\leq\exp\left(-e^{-2r}\int_0^L e^{-2g(t)}dt\right),
\end{align*}
the first term on the right-hand side is bounded by $C(r-x)\Phi_{L;u,v}(r,g,h).$ For the second term, using $1-e^{-y}\leq y$ for $y\geq0$, we obtain
\begin{multline*}
e^{-2(u+v)r-2vg(L)}\left|\exp\left(-e^{-2x}\int_0^L e^{-2g(t)}dt\right)-\exp\left(-e^{-2r}\int_0^L e^{-2g(t)}dt\right)\right|
\\  = e^{-2(u+v)r-2vg(L)}\exp\left(-e^{-2r}\int_0^L e^{-2g(t)}dt\right)\left(1-\exp\left(-\left(e^{-2x}-e^{-2r}\right)\int_0^L e^{-2g(t)}dt\right)\right)
\\ \leq C(r-x)e^{-2(u+v)r-2vg(L)}\exp\left(-e^{-2r}\int_0^L e^{-2g(t)}dt\right)e^{-2r}\int_0^L e^{-2g(t)}dt.
\end{multline*}
Applying $ze^{-z}\leq Ce^{-z/2}$ for $z\geq0$, the last expression is bounded by
\begin{align*}
C(r-x)e^{-2(u+v)r-2vg(L)}\exp\left(-\frac{e^{-2r}}{2}\int_0^L e^{-2g(t)}dt\right) \leq C(r-x)\Phi_{L;u,v}\left(r+2^{-1}\log{(2)},g,h\right).
\end{align*}
Consequently,
\begin{align}\label{e:weightcoordinatebound}
\left|\Phi_{L;u,v}(x,g,h)-\Phi_{L;u,v}(r,g,h) \right| \leq CN^{-1/2}\left(\Phi_{L;u,v}(r,g,h)+\Phi_{L;u,v}\left(r+2^{-1}\log{(2)},g,h\right)\right).
\end{align}

The contribution from the first term in the preceding bound is
\begin{align*}
\mathrm{Lip}(F)\sum_{x\in\mathbb Z^{(N)}}\int_x^{x+N^{-1/2}} (r-x)\mathbb E_{\mathbb W_L}\left[\Phi_{L;u,v}(x,\cdot,\cdot)\right]dr 
&=\frac{\mathrm{Lip}(F)}{2\sqrt N}\left(N^{-1/2}\sum_{x\in\mathbb Z^{(N)}}\mathbb E_{\mathbb W_L}\left[\Phi_{L;u,v}(x,\cdot,\cdot)\right]\right).
\end{align*}
By \eqref{e:continuumlatticesum}, the expression in parentheses is uniformly bounded in $N$. Therefore, this term converges to zero.

By \eqref{e:weightcoordinatebound}, the remaining term is bounded by
\begin{multline*}
C N^{-1/2}\|F\|_{\infty}\sum_{x\in\mathbb Z^{(N)}}\int_x^{x+N^{-1/2}}\mathbb E_{\mathbb W_L}\left[\Phi_{L;u,v}(r,\cdot,\cdot)+\Phi_{L;u,v}\left(r+\frac{\log 2}{2},\cdot,\cdot\right)\right]dr
\\  =C N^{-1/2}\|F\|_{\infty}\int_{\widetilde{x}_N}^{\infty}\mathbb E_{\mathbb W_L}\left[\Phi_{L;u,v}(r,\cdot,\cdot)+\Phi_{L;u,v}\left(r+\frac{\log 2}{2},\cdot,\cdot\right)\right]dr \leq 2C N^{-1/2}\|F\|_{\infty}\mathcal Z_{L;u,v},
\end{multline*}
which converges to zero. We have now shown that all three terms on the right-hand side of \eqref{e:firstbound} converge to zero, which implies \cref{p:almost}.
\end{proof}

\appendix

\section{Technical Lemmas}\label{s:rws}
In this appendix, we prove \cref{l:combapprox} and \cref{l:measurechange}.
\begin{proof}[Proof of \cref{l:combapprox}]
For the duration of the proof, we will write $n:=\lfloor LN\rfloor,d:=\left\lfloor anN^{-1/2}\right\rfloor.$ For all sufficiently large $N$, the hypotheses of \cref{l:combapprox} imply that $|k|<n$ and $|k+d|<n,$ so both binomial coefficients in the proposition are positive. 

If $a=0$, then $d=0$ and the result follows immediately. We may therefore assume that $a\neq0$ and, after increasing $N_{L,a}$ if necessary, that $d\neq0$.

When $d>0$, we have the exact identity
\begin{align}\label{e:binomialratio-positive}
\binom{2n}{n+k}\binom{2n}{n+k+d}^{-1}=\prod_{j=1}^{d}\frac{n+k+j}{n-k-j+1}.
\end{align}
It therefore suffices to show that
\begin{align*}
\exp{\left(-\frac{2ak}{\sqrt N}\right)}\binom{2n}{n+k}\binom{2n}{n+k+d}^{-1}
\end{align*}
is bounded uniformly over the indicated values of $N$ and
$k$. There exists a universal constant $C>0$ such that, for every $|z|\leq1/2$,
\begin{align*}
\frac{1+z}{1-z} \leq \exp{\left(2z+C|z|^3\right)}.
\end{align*}

Suppose first that $d>0$. For $j\in\llbracket1,d\rrbracket$, we have
\begin{align*}
\frac{n+k+j}{n-k-j+1}=\left(1+\dfrac{k+j-1/2}{n+1/2}\right)\left(1-\dfrac{k+j-1/2}{n+1/2}\right)^{-1}.
\end{align*}
Since $|k|<n^{5/6}$ and $|k+d|<n^{5/6}$, we note that $|k+j-1/2|\leq n^{5/6}+1$ for every $j\in\llbracket1,d\rrbracket$. Hence, for all sufficiently large $N$, $|(k+j-1/2)/(n+1/2)|\leq 1/2.$ Using the exact product representation of the binomial ratio, we therefore obtain
\begin{align*}
\binom{2n}{n+k}\binom{2n}{n+k+d}^{-1} & = \prod_{j=1}^{d} \frac{n+k+j}{n-k-j+1}
\\ & \leq\exp{\left(\frac{2}{n+1/2}\sum_{j=1}^{d}\left(k+j-1/2\right)+\frac{C}{(n+1/2)^3}\sum_{j=1}^{d} |k+j-1/2 |^3\right)}
\\ &=\exp{\left(\frac{2dk+d^2}{n+1/2}+\frac{C}{(n+1/2)^3}\sum_{j=1}^{d} |k+j-1/2 |^3\right)}.
\end{align*}

By the definition of $d$, $|d|\leq |a|nN^{-1/2}+1\leq C_{L,a}n^{1/2}$ for all sufficiently large $N$. Consequently, the previous expression is bounded by 
\begin{align*}
\binom{2n}{n+k}\binom{2n}{n+k+d}^{-1}\leq\exp{\left(\frac{2dk+d^2}{n+1/2}+C_{L,a}\right)}.
\end{align*}

We use symmetry of the binomial coefficients to obtain the same estimate when $d<0$.  By the definition of $d$, $|d-anN^{-1/2}|<1.$
Consequently, since $|k|<n^{5/6}$,
\begin{align*}
2|k|\left|\frac{d}{n+1/2}-\frac{a}{\sqrt N}\right| & \leq 2n^{5/6}\left(\frac{1}{n+1/2} +\frac{|a|}{2\sqrt N\left(n+1/2\right)} \right),
\end{align*}
which is uniformly bounded for all sufficiently large $N$. Moreover, $|d|\leq |a|nN^{-1/2}+1,$
and therefore,
\begin{align*}
\frac{d^2}{n+1/2} & \leq \frac{2a^2n^2N^{-1}+2}{n+1/2} \leq 2a^2 nN^{-1}+2 \leq 2a^2L+2.
\end{align*}
Therefore, after enlarging $C_{L,a}$ if necessary, it follows that
\begin{align*}\exp{\left(-\frac{2ak}{\sqrt N}\right)}\binom{2n}{n+k} \binom{2n}{n+k+d}^{-1}\leq e^{C_{L,a}}.
\end{align*}

Recalling that $n=\lfloor LN\rfloor$ and $d=\left\lfloor a\lfloor LN\rfloor N^{-1/2}\right\rfloor,$ and setting $D_{L,a}:=e^{C_{L,a}}$, we obtain
\begin{align*}
\exp{\left(-\frac{2ak}{\sqrt N}\right)}\binom{2\lfloor LN\rfloor}{\lfloor LN\rfloor+k} \leq D_{L,a}\binom{2\lfloor LN\rfloor}{
\lfloor LN\rfloor+k+
\left\lfloor
a\lfloor LN\rfloor N^{-1/2}
\right\rfloor}.
\end{align*}
This proves \cref{l:combapprox}.
\end{proof}
Finally, we give the proof of \cref{l:measurechange}. 
\begin{proof}[Proof of \cref{l:measurechange}] Let $\mathbb{W}_{L,m}$ denote the law of a Brownian motion with
drift $m$ and diffusion coefficient $1/\sqrt{2}$. By the
Cameron--Martin formula,
\begin{align*}
Z_t(g) := \frac{d\mathbb{W}_{L,m}|_{\mathcal{F}_t}}{d\mathbb{W}_L|_{\mathcal F_t}}(g) =\exp{\left( 2mg(t)-m^2t \right)}.
\end{align*}
Under $\mathbb{W}_{L,m}$, the translated path $g_{m}(t)=g(t)-mt$ has law $\mathbb{W}_L$. Since $g(L)$ is Gaussian with mean $0$ and variance $L/2$, for $r\in\mathbb{R}$ its density is $\rho_{g(L)}(r)=(\pi L)^{-1/2}\exp{\left(-r^{2}L^{-1}\right)}.$ Therefore, for every Borel set $A$,
\begin{align}\label{e:bd}
    \mathbb{E}_{\mathbb{W}_{L}}\left[F(g)\mathbf{1}_{\{g(L)\in A\}}\right]=\int_{A}\mathbb{E}_{\mathbb{W}_{L}}\left[F(g)\big|g(L)=r\right]\rho_{g(L)}(r)dr.
\end{align}
Here, the conditional expectation is understood with respect to the Brownian bridge law with endpoint $r$. We now perform the change-of-measure calculation. Let $A$ be a Borel set, and write $A+mL:=\{a+mL:a\in A\}.$ Since $g_m$ under $\mathbb{W}_{L,m}$ has law $\mathbb{W}_L$, and since $g_m(L)=g(L)-mL$, we have
\begin{align*} 
    \mathbb{E}_{\mathbb{W}_L}\left[F(g)\mathbf 1_{\{g(L)\in A\}}\right] &=\mathbb{E}_{\mathbb{W}_{L,m}}\left[F(g_m)\mathbf 1_{\{g_m(L)\in A\}}\right]
    \\ & =\mathbb{E}_{\mathbb{W}_{L,m}}\left[F(g_m)\mathbf 1_{\{g(L)\in A+mL\}}\right]
    \\ &= \mathbb{E}_{\mathbb{W}_L}\left[F(g_m)Z_L(g)\mathbf 1_{\{g(L)\in A+mL\}}\right].
\end{align*} 
Applying \eqref{e:bd} to the first expression, we obtain
\begin{align*}
\mathbb{E}_{\mathbb{W}_L}\left[ F(g)\mathbf 1_{\{g(L)\in A\}} \right] = \int_A \mathbb{E}_{\mathbb{W}_L}\left[F(g)\big|g(L)=a\right]\rho_{g(L)}(a) da.
\end{align*}
On the other hand, since $Z_L(g)=\exp{\left(2mg(L)-m^2L\right)},$ the preceding calculation yields 
\begin{align*}
    \mathbb{E}_{\mathbb{W}_L}\left[F(g_m)Z_L(g)\mathbf 1_{\{g(L)\in A+mL\}}\right]
    & =\int_{A+mL}e^{2mr-m^2L}\mathbb{E}_{\mathbb{W}_L}\left[F(g_m)\big |g(L)=r\right]\rho_{g(L)}(r)dr.
\end{align*}
After making the change of variables $r=a+mL$, this becomes
\begin{align*}
    \int_A e^{2m(a+mL)-m^2L}\mathbb{E}_{\mathbb{W}_L}\left[F(g_m)\big |g(L)=a+mL\right] \rho_{g(L)}(a+mL) da.
\end{align*}
Since
\begin{align*}
e^{2m(a+mL)-m^2L}\rho_{g(L)}(a+mL)=\rho_{g(L)}(a),
\end{align*} we conclude that for every Borel set $A$,
\begin{align*}
    \int_A\mathbb{E}_{\mathbb{W}_L}\left[F(g)\big|g(L)=a\right]\rho_{g(L)}(a)da & =\int_A\mathbb{E}_{\mathbb{W}_L}\left[F(g_m)\big|g(L)=a+mL\right]\rho_{g(L)}(a)da.
\end{align*}
It follows that, for Lebesgue-almost every $a\in\mathbb R$,
\begin{align*}
\mathbb{E}_{\mathbb{W}_L}\left[F(g)\big|g(L)=a\right]=\mathbb{E}_{\mathbb{W}_L}\left[F(g_m)\big|g(L)=a+mL\right].
\end{align*}
Both sides are continuous functions of $a$. This proves the claim.
\end{proof}
\section{Pointwise Convergence}\label{a:pointwise}
In this appendix, we follow the argument in \cite{BD23} to prove the pointwise convergence of the unnormalized weights. We quote a proposition from Durrett~\cite{Dur19}.
\begin{proposition}[Exercise 3.1.1~\cite{Dur19}]\label{p:dur}
    If the following three conditions are met
    \begin{align*}\lim_{n\to\infty}\max_{0\leq j\leq n}|c_{j,n}|=0, & & \lim_{n\to\infty}\sum_{j=0}^{n}c_{j,n}=\lambda, & & \sup_{n}\sum_{j=0}^{n}|c_{j,n}|<\infty,\end{align*} with $c_{j,n}\in\mathbb{R}$ for every $n\in\mathbb{Z}_{>0}$ and $j\in\llbracket0,n\rrbracket$, then $\lim_{n\to\infty}\prod_{j=0}^{n}(1+c_{j,n})=e^{\lambda}.$
\end{proposition}
We use the immediate triangular-array variant of \cref{p:dur} in which the $N$-th row is indexed by $i\in\llbracket0,\lfloor LN\rfloor\rrbracket.$ The proof is unchanged.
\begin{proposition}
Fix $L>0$ and $u,v\in\mathbb{R}$. Then, as $N\to\infty$, for every $(x,g,h)\in\mathbb{R}\times C_0[0,L]\times C_0[0,L],$ $\Phi^{(N)}_{L;u,v}(x,g,h)\longrightarrow\Phi_{L;u,v}(x,g,h)$.
\end{proposition}
\begin{proof}
For the purposes of this appendix, for every $N$ such that $\lfloor LN\rfloor\geq1$, we extend $\Phi^{(N)}_{L;u,v}$ to $\mathbb{R}\times C_0[0,L]\times C_0[0,L]$ using the unrestricted product formula
\begin{align*}
\Phi^{(N)}_{L;u,v}(x,g,h) &= \exp{\left(-2(u+v)x-2vg(L)\right)} \prod_{i=0}^{\lfloor LN\rfloor} \left(1-\frac{e^{-2(g(t_i)+x)}}{N}\right),
\\ \Phi_{L;u,v}(x,g,h) & = \exp{\left(-2(u+v)x-2vg(L)-e^{-2x}\int_0^L e^{-2g(t)} dt\right)}.
\end{align*}
We consider the sequence $c_{i,N}:=-N^{-1}e^{-2(g(t_i)+x)}$ for $i\in\llbracket0,\lfloor LN\rfloor\rrbracket$.
    For any fixed function $g$ on the interval $[0,L]$, the function $e^{-2(g(t)+x)}$ is bounded by a constant $C_{g,x}>0$. Therefore,  $$\lim_{N\to\infty}\max_{0\leq i\leq\lfloor LN\rfloor}|c_{i,N}|\leq\lim_{N\to\infty}N^{-1}C_g=0.$$
Likewise, because $g\in C_{0}[0,L]$, we also know that the function $t\mapsto e^{-2(g(t)+x)}$ belongs to $C[0,L]$. Therefore, by the definition of the Riemann integral,
\begin{align*}
\lim_{N\to\infty}\sum_{i=0}^{\lfloor LN\rfloor}\frac{e^{-2(g(t_i)+x)}}{N} & =\lim_{N\to\infty}\left(\frac{\lfloor LN\rfloor}{LN}\left(\frac{L}{\lfloor LN\rfloor}\sum_{i=0}^{\lfloor LN\rfloor-1}e^{-2(g(t_i)+x)}\right)+\frac{e^{-2(g(L)+x)}}{N}\right)=e^{-2x}\int_{0}^{L}e^{-2g(t)}dt.
\end{align*}
Finally, by the same reasoning,
\begin{align*}
\lim_{N\to\infty}\sum_{i=0}^{\lfloor LN\rfloor}|c_{i,N}| & = \lim_{N\to\infty}\sum_{i=0}^{\lfloor LN\rfloor} N^{-1}e^{-2(g(t_i)+x)}=e^{-2x}\int_{0}^{L}e^{-2g(t)}dt,
\end{align*}
and therefore $\sup_{N: \lfloor LN\rfloor\geq1}\sum_{i=0}^{\lfloor LN\rfloor}|c_{i,N}|<\infty.$
    Putting all of these calculations together, \cref{p:dur} implies
    \begin{align*}
        \lim_{N\to\infty} e^{-2(u+v)x-2vg(L)}\prod_{i=0}^{\lfloor LN\rfloor}\left(1-\frac{e^{-2(g(t_{i})+x)}}{N}\right) 
         & = \exp{\left(-2(u+v)x-2vg(L)-e^{-2x}\int_{0}^{L}e^{-2g(t)}dt\right)}.
    \end{align*}
\end{proof}

\bibliographystyle{alpha}
\bibliography{mybibliography.bib}

\begin{thebibliography}{BKWW23}

\bibitem[Bil99]{Bill99}
Patrick Billingsley.
\newblock {\em Convergence of probability measures}.
\newblock Wiley Series in Probability and Statistics: Probability and Statistics. John Wiley \& Sons, Inc., New York, second edition, 1999.
\newblock A Wiley-Interscience Publication.

\bibitem[BKWW23]{BKWW23}
W{\l}odek Bryc, Alexey Kuznetsov, Yizao Wang, and Jacek Weso{\l}owski.
\newblock Markov processes related to the stationary measure for the open {KPZ} equation.
\newblock {\em Probab. Theory Related Fields}, 185(1-2):353--389, 2023.

\bibitem[BLD22]{BLD22}
Guillaume Barraquand and Pierre Le~Doussal.
\newblock Steady state of the {KPZ} equation on an interval and {L}iouville quantum mechanics.
\newblock {\em Europhysics Letters}, 137(6):61003, March 2022.

\bibitem[BLD23]{BD23}
Guillaume Barraquand and Pierre Le~Doussal.
\newblock Stationary measures of the {KPZ} equation on an interval from {E}naud-{D}errida's matrix product ansatz representation.
\newblock {\em J. Phys. A}, 56(14):Paper No. 144003, 14, 2023.

\bibitem[CK24]{CK21}
Ivan Corwin and Alisa Knizel.
\newblock Stationary measure for the open {KPZ} equation.
\newblock {\em Comm. Pure Appl. Math.}, 77(4):2183--2267, 2024.

\bibitem[CN20]{CN20}
Sylvie Corteel and Arthur Nunge.
\newblock Combinatorics of the 2-species exclusion processes, marked {L}aguerre histories, and partially signed permutations.
\newblock {\em Electron. J. Combin.}, 27(2):Paper No. 2.53, 27, 2020.

\bibitem[Cor22]{C22}
Ivan Corwin.
\newblock Some recent progress on the stationary measure for the open {KPZ} equation.
\newblock In {\em Toeplitz operators and random matrices---in memory of {H}arold {W}idom}, volume 289 of {\em Oper. Theory Adv. Appl.}, pages 321--360. Birkh\"{a}user/Springer, Cham, 2022.

\bibitem[CS18]{CS18}
Ivan Corwin and Hao Shen.
\newblock Open {ASEP} in the weakly asymmetric regime.
\newblock {\em Comm. Pure Appl. Math.}, 71(10):2065--2128, 2018.

\bibitem[DEHP93]{DEHP93}
B.~Derrida, M.~R. Evans, V.~Hakim, and V.~Pasquier.
\newblock Exact solution of a 1{D} asymmetric exclusion model using a matrix formulation.
\newblock {\em J. Phys. A}, 26(7):1493, 1993.

\bibitem[DEL04]{DEL}
B.~Derrida, C.~Enaud, and J.~L. Lebowitz.
\newblock The asymmetric exclusion process and {B}rownian excursions.
\newblock {\em J. Statist. Phys.}, 115(1-2):365--382, 2004.

\bibitem[Dur19]{Dur19}
Rick Durrett.
\newblock {\em Probability---theory and examples}, volume~49 of {\em Cambridge Series in Statistical and Probabilistic Mathematics}.
\newblock Cambridge University Press, Cambridge, fifth edition, 2019.

\bibitem[DW21]{DW}
Evgeni Dimitrov and Xuan Wu.
\newblock K{MT} coupling for random walk bridges.
\newblock {\em Probab. Theory Related Fields}, 179(3-4):649--732, 2021.

\bibitem[ED04]{ED04}
C.~Enaud and B.~Derrida.
\newblock Large deviation functional of the weakly asymmetric exclusion process.
\newblock {\em J. Statist. Phys.}, 114(3-4):537--562, 2004.

\bibitem[ER96]{ER96}
Fabian H.~L. Essler and Vladimir Rittenberg.
\newblock Representations of the quadratic algebra and partially asymmetric diffusion with open boundaries.
\newblock {\em J. Phys. A}, 29(13):3375--3407, 1996.

\bibitem[GR15]{GR}
I.~S. Gradshteyn and I.~M. Ryzhik.
\newblock {\em Table of integrals, series, and products}.
\newblock Elsevier/Academic Press, Amsterdam, eighth edition, 2015.
\newblock Translated from the Russian, Translation edited and with a preface by Daniel Zwillinger and Victor Moll.

\bibitem[KM22]{KM22}
Alisa Knizel and Konstantin Matetski.
\newblock The strong {F}eller property of the open {KPZ} equation.
\newblock {\em arXiv:2211.04466}, 2022.

\bibitem[KMT75]{KMT75}
J.~Koml\'{o}s, P.~Major, and G.~Tusn\'{a}dy.
\newblock An approximation of partial sums of independent {${\rm RV}$}'s and the sample {${\rm DF}$}. {I}.
\newblock {\em Z. Wahrscheinlichkeitstheorie und Verw. Gebiete}, 32:111--131, 1975.

\bibitem[KPZ86]{KPZ86}
Mehran Kardar, Giorgio Parisi, and Yi-Cheng Zhang.
\newblock Dynamic scaling of growing interfaces.
\newblock {\em Phys. Rev. Lett.}, 56:889--892, Mar 1986.

\bibitem[LL10]{L10}
Gregory~F. Lawler and Vlada Limic.
\newblock {\em Random walk: a modern introduction}, volume 123 of {\em Cambridge Studies in Advanced Mathematics}.
\newblock Cambridge University Press, Cambridge, 2010.

\bibitem[MS97]{MS97}
K.~Mallick and S.~Sandow.
\newblock Finite-dimensional representations of the quadratic algebra: applications to the exclusion process.
\newblock {\em J. Phys. A}, 30(13):4513--4526, 1997.

\bibitem[Mue91]{M91}
Carl Mueller.
\newblock On the support of solutions to the heat equation with noise.
\newblock {\em Stochastics Stochastics Rep.}, 37(4):225--245, 1991.

\bibitem[MY05]{MY05}
Hiroyuki Matsumoto and Marc Yor.
\newblock Exponential functionals of {B}rownian motion. {I}. {P}robability laws at fixed time.
\newblock {\em Probab. Surv.}, 2:312--347, 2005.

\bibitem[Par19]{P19}
Shalin Parekh.
\newblock The {KPZ} limit of {ASEP} with boundary.
\newblock {\em Comm. Math. Phys.}, 365(2):569--649, 2019.

\bibitem[Par22]{P22}
Shalin Parekh.
\newblock Ergodicity results for the open {KPZ} equation.
\newblock {\em arXiv:2212.08248}, 2022.

\bibitem[Sas99]{Sa99}
Tomohiro Sasamoto.
\newblock One-dimensional partially asymmetric simple exclusion process with open boundaries: orthogonal polynomials approach.
\newblock {\em J. Phys. A}, 32(41):7109--7131, 1999.

\bibitem[USW04]{USW}
Masaru Uchiyama, Tomohiro Sasamoto, and Miki Wadati.
\newblock Asymmetric simple exclusion process with open boundaries and {A}skey-{W}ilson polynomials.
\newblock {\em J. Phys. A}, 37(18):4985--5002, 2004.

\bibitem[Yor92]{Y92}
Marc Yor.
\newblock On some exponential functionals of {B}rownian motion.
\newblock {\em Adv. in Appl. Probab.}, 24(3):509--531, 1992.

\end{thebibliography}
\end{document}